\documentclass[11pt]{amsart}
\usepackage[utf8]{inputenc}
\usepackage{setspace}
\usepackage[margin=1in]{geometry}
\usepackage{cite, comment, url}
\usepackage{tikz, graphicx, marvosym, subfigure, float, circuitikz}
\usetikzlibrary{decorations.pathreplacing}
\usepackage{dsfont, amsfonts, amssymb, amsmath, amsthm, mathtools}
\usepackage{algorithm, algpseudocode}

\usepackage[breaklinks,backref=page]{hyperref}
\usepackage{cleveref}
\hypersetup{
	colorlinks = true, 
	urlcolor = cyan, 
	linkcolor = teal, 
	citecolor = cyan 
}

\renewcommand*{\backref}[1]{}
\renewcommand*{\backrefalt}[4]{
	\ifcase #1 \textcolor{red}{Not cited.}%
	\or $\uparrow$#2%
	\else $\uparrow$#2%
	\fi%
}

\makeatletter
\renewcommand{\l@section}{%
  \@tocline{1}{0em}{0em}{}{}%
}
\renewcommand{\l@subsection}{%
  \@tocline{2}{-0.5\baselineskip}{1.7em}{}{}%
}
\makeatother
\numberwithin{equation}{section}

\newtheorem{theorem}{Theorem}[section]

\newtheorem{lemma}[theorem]{Lemma}
\newtheorem{proposition}[theorem]{Proposition}

\newtheorem{question}{Question}

\renewcommand{\pmod}[1]{\ (\mathrm{mod}\ #1)}

\newcommand{\sm}{\!\setminus\!}
\newcommand{\R}{\mathbb R}

\newcommand{\N}{\mathbb N}
\newcommand{\cH}{\mathcal H}
\newcommand{\cF}{\mathcal F}
\newcommand{\cV}{\mathcal V}

\newcommand{\cQ}{\mathcal Q}
\newcommand{\x}{\mathbf x}
\newcommand{\y}{\mathbf y}
\newcommand{\z}{\mathbf z}
\newcommand{\uu}{\mathbf u}
\newcommand{\vv}{\mathbf v}

\newcommand{\aaa}{\mathbf a}
\newcommand{\ex}{{\rm ex}}

\newcommand{\ind}[1]{\smash{\xrightarrow[\smash{\raisebox{0.7ex}{\tiny $ind$}}]{#1}}}

\newcommand{\chigraph}[2]{\chi_{#2}(#1)}
\newcommand{\chigraphind}[2]{\chi^{\mathrm{ind}}_{#2}(#1)}
\newcommand{\alphagraphind}[2]{\alpha^{\mathrm{ind}}_{#2}(#1)}

\title{Ramsey problems for graphs in Euclidean spaces and Cartesian powers}

\author{Maria Axenovich}

\author{Dingyuan Liu}

\author{Arsenii Sagdeev}

\address{Karlsruhe Institute of Technology, Englerstraße 2, D-76131 Karlsruhe, Germany}
\email{maria.aksenovich@kit.edu, liu@mathe.berlin, sagdeevarsenii@gmail.com}

\begin{document}

\begin{abstract}
Given a graph $H$, let $\chi_H(\mathbb{R}^n)$ be the smallest positive integer $r$ such that there exists an $r$-coloring of $\mathbb{R}^n$ with no monochromatic unit-copy of $H$, that is a set of $|V(H)|$ vertices of the same color such that any two vertices corresponding to an edge of $H$ are at distance one.

This Ramsey-type function extends the famous Hadwiger--Nelson problem on the chromatic number $\chi(\mathbb{R}^n)=\chi_{K_2}(\mathbb{R}^n)$ of the space from a complete graph $K_2$ on two vertices to an arbitrary graph $H$. It also extends the classical Euclidean Ramsey problem for congruent monochromatic subsets to the family of those defined by a specific subset of unit distances.

Among others, we show that $\chi_H(\mathbb{R}^n)=\chi(\mathbb{R}^n)$ for any even cycle $H$ of length $8$ or at least $12$ as well as for any forest and that $\chi_H(\mathbb{R}^n)=\lceil\chi(\mathbb{R}^n)/2\rceil$ for any sufficiently long odd cycle. Our main tools and results, which are of independent interest, establish that Cartesian powers enjoy Ramsey-type properties for graphs with favorable Tur\'an-type characteristics, such as zero hypercube Tur\'an density.

In addition, we prove induced variants of these results, find bounds on $\chi_H(\mathbb{R}^n)$ for growing dimensions $n$, and prove a canonical-type result. We conclude with many open problems. One of these is to determine $\chi_{C_4}(\mathbb{R}^2)$, for a cycle $C_4$ on four vertices.
    
    {\centering
    \includegraphics[scale=0.7,width=0.93\linewidth]{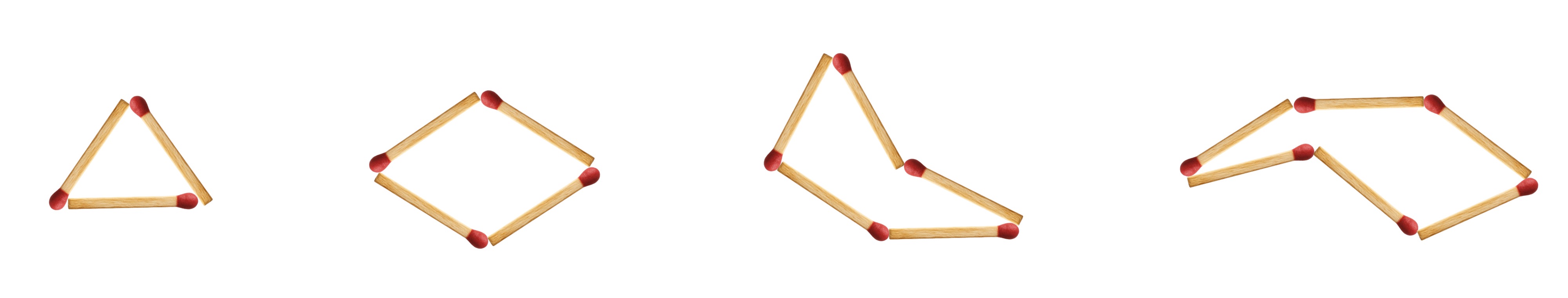}}
\end{abstract}

\maketitle

\vspace{-3em}
\tableofcontents

\section{Introduction} \label{sec:1}

\subsection{The chromatic number $\chigraph{\R^n}{H}$ for a graph $H$, its values for forests and cycles} \label{sec:1.1}

Ramsey theory, see Graham, Rothschild, and Spencer~\cite{GRS91}, is an esteemed branch of modern combinatorics that has been extensively studied over the last century. An important class of geometric Ramsey problems originates from the classic papers~\cite{EGMRSS1, EGMRSS2, EGMRSS3} by Erd\H{o}s, Graham, Montgomery, Rothschild, Spencer, and Straus who defined a set $A \subseteq \R^d$ to be {\it Ramsey} if $\chi(\R^n,A)$ tends to infinity with $n$, where $\chi(\R^n,A)$ is defined as the smallest $r\in \N$ such that there exists an $r$-coloring of $\R^n$ with no monochromatic isometric copy of $A$.
When $A$ is a set of two points at a unit distance apart, $\chi(\R^n,A)$ is also known as the \emph{chromatic number of $\R^n$}, abbreviated as $\chi(\R^n)$.
The inequalities $\chi(\R^n,A)>r$ and $\chi(\R^n,A)\le r$ are often denoted by $\R^n \xrightarrow{r} A$ and $\R^n \not\xrightarrow{r} A$, respectively. 
The main ``meta-question" of Euclidean Ramsey theory is to classify Ramsey sets, see also  Section \ref{sec:background}.

We study a different problem of a similar nature.
For a graph $H=(V, E)$, $E\subseteq \binom{V}{2}$, we say that a set $A\subseteq \R^n$ with $|A|=|V|$ {\it represents} $H$ or is a {\it unit-copy} of $H$ if 
there is a bijection $f:V\rightarrow A$ such that $uv\in E$ implies that $||f(u)- f(v)||=1$, for the Euclidean norm $||\cdot||$. Such a unit-copy is {\it induced} if for any $uv\not\in E$, $||f(u)- f(v)||\neq 1$.
We call graphs that have an induced unit-copy in $\R^n$ {\it unit-distance graphs in $\R^n$}. 
Note that each finite graph has an induced unit-copy in $\R^n$ for some $n$, see for example Alon and Kupavskii~\cite{AK14} and Frankl,  Kupavskii, and  Swanepoel~\cite{NKS}.

Throughout this paper, graphs are assumed to be non-empty and finite unless explicitly stated. We denote a clique on $d$ vertices by $K_d$ and a cycle of length $\ell$ by $C_\ell
$.
Given a graph $H$, let $\chigraph{\R^n}{H}$ be the smallest $r\in \N$ such that there exists an $r$-coloring of $\R^n$ with no monochromatic unit-copy of $H$.
The function $\chigraphind{\R^n}{H}$ is defined analogously for induced unit-copy of $H$. 
We use the classical arrow-expressions  $\R^n \xrightarrow{r} H$ and $\R^n \not\xrightarrow{r} H$ that are equivalent to the statements $\chigraph{\R^n}{H}>r$ and $\chigraph{\R^n}{H}\le r$, respectively.  
If $H=K_2$, then $\chigraphind{\R^n}{H} = \chigraph{\R^n}{H} = \chi(\R^n)$ by definition. 

We are looking for monochromatic sets of vertices that preserve unit distances between some vertices, but not necessarily all of them. One can also visualize this setting by considering the ``match-stick" model of a given unit-distance graph, where the edges are rigid ``sticks" with flexible ``joints" at the vertices.  Our goal is to find some monochromatic copy of that model.
For example, if $\triangle$ is a set  of $3$ points that are pairwise at distance $1$ and $H$ is a complete graph on $3$ vertices, then $\chi(\R^2,\triangle)=\chigraph{\R^2}{H}$, since each unit-copy of $K_3$ is a copy of $\triangle$. In this case our setting matches the classical Euclidean Ramsey one. Moreover, a red-blue coloring of $\R^2$ formed by iterated red and blue infinite stripes of width $\sqrt{3}/2$ has no monochromatic copy of $\triangle$, thus giving that $\chigraph{\R^2}{H}=2$. On the other hand, when $H$ is a path on $3$ vertices, we have a large family of $3$-point configurations in the plane that are unit-copies of $H$. One can easily show that $\chigraph{\R^2}{H} > 2$ because for any red-blue coloring of $\R^2$ and any, say, red vertex $v$, a unit circle centered at $v$ either has two red vertices, thus giving us a red unit-copy of $H$, or all but one blue vertices, giving a blue unit-copy of $H$. The first  result of this paper  concerns forests and cycles.

\begin{theorem}
\label{thm:main}
Let $n\geq 2$ be an integer. Then the following holds:
\begin{enumerate}
    \setlength\itemsep{1ex}   
    \item $\chigraphind{\R^n}{F} = \chigraph{\R^n}{F} =\chi(\R^n)$ for any forest $F$.
    \item $\chigraphind{\R^n}{C_{2\ell}} = \chigraph{\R^n}{C_{2\ell}} =\chi(\R^n)$ for $\ell=4$ and any integer $\ell \geq 6$.
    \item There exists $\ell_0=\ell_0(n)$ such that $\chigraph{\R^n}{C_{2\ell+1}} =  \lceil \chi(\R^n)/ 2 \rceil$ for any integer $\ell \geq \ell_0$.
\end{enumerate}
\end{theorem}

\Cref{thm:main} is a consequence of  more general results that we  present in ~\Cref{sec:1.2,sec:1.3}.  We give a background on Euclidean Ramsey theory and make some basic observations on the functions $\chigraph{\R^n}{H}$ and $\chigraphind{\R^n}{H}$ in Sections~\ref{sec:background} and~\ref{sec:basic-observations}, respectively. We conclude the introductory chapter with Section \ref{sec:growing} stating results on growing dimensions and a canonical-type result.

\subsection{Ramsey and Tur\'an properties of large Cartesian powers of graphs} \label{sec:1.2}

The main tool we use to prove \Cref{thm:main} is a certain Ramsey-type result on Cartesian powers of graphs that may be of independent interest. We need several definitions to state this general result.

For $N \in \N$, let $[N]\coloneq\{1,\dots,N\}$. Given two graphs $G$ and $H$, their \emph{Cartesian product $G\square H$} is defined as the graph on the vertex set $V(G)\times V(H)$, where the vertices $(x_1,x_2),(y_1,y_2)$ form an edge in $G\square H$ if and only if either $x_1y_1\in E(G)$ and $x_2=y_2$, or $x_1=y_1$ and $x_2y_2\in E(H)$. Let $G^{\square 1}=G$ and for $N\geq 2$, let $G^{\square N}= G^{\square N-1} \square G$ be the $N$th {\it Cartesian power} of $G$.

A \emph{$G$-slice} of $G^{\square N}$ is a subgraph induced by the Cartesian product of $V(G)$ and $N-1$ one-element subsets of $V(G)$ in some order. There are exactly $N|V(G)|^{N-1}$ $G$-slices in $G^{\square N}$, each of which is isomorphic to $G$.
We say that a graph $H$ is of \emph{zero (induced) $G$-slice density} if for every $\varepsilon>0$, there exists $N_0=N_0(\varepsilon)$ such that for each integer $N\ge N_0$, every subgraph of $G^{\square N}$ that contains at least an $\varepsilon$-fraction of its $G$-slices contains an (induced) copy of $H$. Given a graph $G$, a collection $\mathcal{F}$ of graphs, and a positive integer $r$, we write $G\xrightarrow{r}\mathcal{F}$ (resp. $G \ind{r}\mathcal{F}$) if every $r$-coloring of $V(G)$ contains a monochromatic copy (resp. induced copy) of some graph in $\mathcal{F}$.

\begin{theorem} \label{thm:Turan}    
Let $r$ be a positive integer, $G_1,\dots,G_m$, $H_1,\dots,H_m$, and $G$ be graphs.
If  $H_i$ has zero  $G_i$-slice density for each $i \in [m]$ and $G\xrightarrow{r}\{G_1,\dots,G_m\}$, then   $G^{\square N} \xrightarrow{r}\{H_1,\dots,H_m\}$ for a sufficiently large integer~$N$.
Similarly, if $H_i$ has zero induced  $G_i$-slice density for each $i \in [m]$ and  $G \ind{r}\{G_1,\dots,G_m\}$, then $G^{\square N} \ind{r}\{H_1,\dots,H_m\}$ for a sufficiently large integer $N$. 
\end{theorem}

In the next result we shall provide several concrete instances of~\Cref{thm:Turan} when $m=1$ by exploring special classes of graphs of zero $G$-slice density. For this we need some definitions.

A \emph{hypercube} $Q_N$ is a graph isomorphic to $K_2^{\square N}$ on the vertex set $\{0,1\}^N$ where two vertices are adjacent if and only if they differ in exactly one coordinate. The \emph{$k$th vertex layer} of $Q_N$ is the set of vertices containing exactly $k$ ones,  the \emph{$k$th edge layer} of $Q_N$ is the subgraph of $Q_N$ induced by the $k$th and $(k-1)$th vertex layers.  A graph $H$ is of \emph{zero (induced) hypercube Tur\'an density} if the maximum fraction of edges in a subgraph of $Q_N$ that does not contain an (induced) copy of $H$ tends to 0 as $N \to \infty$. Conlon \cite{C} described a wide class of graphs that have zero hypercube Tur\'an density,  see also the papers by Zhu and the first author~\cite{Z, A} for more such graphs. Observe that the notion of a zero $G$-slice density generalizes the classical zero hypercube Tur\'an density notion when $G=K_2$.

For a bipartite graph $\Gamma$ with parts $A$ an $B$, define $m(N,\Gamma,A,B)$ to be the maximum number of edges in a subgraph of $Q_N$ that contains no copy of $\Gamma$ with the vertices corresponding to $A$ in the $k$th vertex layer and vertices corresponding to $B$ in $(k-1)$th vertex layer, for any $k\leq N$. We say that $\Gamma$ is of \emph{zero strong hypercube Tur\'an density} with respect to $(A,B)$, if $m(N,\Gamma,A,B) = o(|E(Q_N)|)$ as $N \to \infty$.
Note that it was observed by the first author in~\cite{A} that almost all known graphs of zero hypercube Tur\'an density (i.e., those described by Conlon~\cite{C}) have zero strong hypercube Tur\'an density with respect to any of their bipartition.

Let $H$ be a graph and $u,v\in V(H)$ be any two vertices. Given a bipartite graph $\Gamma$ with parts $A$ and $B$, we define $\Gamma_{A,B}(H,u,v)$ to be the graph obtained from $\Gamma$ by replacing every edge $ab\in E(\Gamma)$, $a\in A$ and $b\in B$, with a copy of $H$, where $(a,b)$ corresponds to $(u,v)$.
A graph $F$ is called an \emph{$H$-forest} if $F$ is a union of pairwise edge-disjoint copies $H_1,\dots,H_m$ of $H$ such that every $V(H_i)$ intersects $\bigcup_{j<i}V(H_j)$ in at most one vertex. Note that $K_2$-forests are simply forests.

\begin{proposition} \label{pro:Turan}
Let $G, H$ be graphs and  $r$ be a positive integer. Then the following holds:
\begin{enumerate}
    \setlength\itemsep{1ex}
    \item If $H$ is vertex-transitive, $F$ is an $H$-forest, and 
    $G \xrightarrow{r}  H$ (resp., $G \ind{r} H$),
     then there is a sufficiently large $N$ such that $G^{\square N} \xrightarrow{r} F$ (resp., $G^{\square N} \ind{r} F$).

    \item If $H$ has zero (induced) hypercube Tur\'an density and $\chi(G)=r+1$, then there is a sufficiently large $N$ such that  $G^{\square N} \xrightarrow{r} H$ (resp., 
    $G^{\square N} \ind{r} H$).
     
    \item Let $u,v$ be two vertices of $H$. Let $\Gamma$ be a bipartite graph with parts $A$ and $B$ of zero strong hypercube Tur\'an density with respect to $(A,B)$. 
   If $G \xrightarrow{r}  H$, then  there is a sufficiently large $N$ such that $G^{\square N} \xrightarrow{r} \Gamma_{A,B}(H,u,v)$.
\end{enumerate}
\end{proposition}

 \Cref{thm:Turan} provides results for large classes of graphs, in particular 
for most even cycles.  An even cycle $C_{2\ell}$ is of zero induced hypercube Tur\'an density if $\ell=4$ or $\ell\ge 6$, which is already implicit from the proof by Chung \cite{Chung92} for even $\ell$. For odd $\ell$ the non-induced version is proved by F\"uredi and Ozkahya \cite{FO}. We provide a proof for the induced version of these statements  in \Cref{app:hypercube-Turan}. Note that $C_4$, $C_6$, and $C_{10}$ are not of zero hypercube Tur\'an density, as follows from constructions by Chung, Conder, Grebennikov and Marciano, see ~\cite{Chung92,conder,GM25}, respectively.

However,~\Cref{thm:Turan} does not cover another important class of so-called \emph{(induced) layered graphs}, that are (induced) subgraphs of an edge layer of a hypercube. For some results on layered graphs, see~\cite{AMW, BLMW}. A characterization of these graphs in terms of edge-colorings is given by Havel and Mor\'avek~\cite{HM}.
Although every graph of zero hypercube Tur\'an density is layered, the other direction is not necessarily true, see for example $C_6$ and $C_{10}$. Our next result shows that monochromatic copies of layered graphs are unavoidable in $2$-colorings of the large Cartesian power of $K_3$.

\begin{theorem}
\label{lem:layered}
    If $H$ is an (induced) layered graph, then there is a sufficiently large $N$ such that $K_3^{\square N} \xrightarrow{2} H$ (resp., $K_3^{\square N} \ind{2} H$).
\end{theorem}

A result of Horvat and Pisanski~\cite[Theorem~3.4]{HP10}, stating that the Cartesian product of unit-distance graphs in $\R^n$ is again a unit-distance graph in $\R^n$, allows us to use our Ramsey-type results for Cartesian powers in Euclidean Ramsey theory. We summarize these results next.

\subsection{General results on $\chigraph{\R^n}{H}$ and $\chigraphind{\R^n}{H}$ for other classes of graphs} \label{sec:1.3}
The following result is almost an immediate consequence of \Cref{pro:Turan}.
\begin{theorem}
\label{general}
Let $H$ be a graph and $n\geq2$ be an integer. Then the following holds:
\begin{enumerate}
    \setlength\itemsep{1ex}
    \item If $H$ is vertex transitive and $F$ is an $H$-forest, then $\chigraph{\R^n}{F}=\chigraph{\R^n}{H}$ and $\chigraphind{\R^n}{F}=\chigraphind{\R^n}{H}$.
    \item If $H$ has zero (induced) hypercube Tur\'an density, then $\chigraph{\R^n}{H}=\chigraph{\R^n}{}$ (resp. $\chigraphind{\R^n}{H}=\chigraph{\R^n}{}$).
    \item Let $u,v$ be two vertices of $H$. Let $\Gamma$ be a bipartite graph with parts $A$ and $B$ of zero strong hypercube Tur\'an density with respect to $(A,B)$ and let $F=\Gamma_{A,B}(H,u,v)$. Then $\chigraph{\R^n}{F}=\chigraph{\R^n}{H}$.
\end{enumerate}
\end{theorem}

Observe that the equality $\chi(\R^2,\triangle)=2$ discussed in \Cref{sec:1.1} implies that $\chigraph{\R^2}{H}=2$ for any unit-distance graph $H$ in the plane that contains a triangle. Our next result shows that $\chigraph{\R^2}{H}=2$ holds for some bipartite graphs $H$ as well, while for many others, including $C_6$ and $C_{10}$, \Cref{lem:layered} implies that $\chigraph{\R^2}{H}>2$.

\begin{theorem}
\label{thm:main-layer}
\
\begin{enumerate}
    \setlength\itemsep{1ex}   
    \item If $H$ is an (induced) layered graph, then $\chigraph{\R^2}{H}>2$ (resp. $\chigraphind{\R^2}{H}>2$). 
    \item There exists a bipartite graph $H$, namely $H=Q_{11}$, such that $\chigraph{\R^2}{H}=2$.
\end{enumerate}
\end{theorem}

The next result addresses the smallest bipartite graph not covered so far, that is the four-cycle. 

\begin{proposition}
\label{leftovers}
We have the following bounds for $C_4$: \textup{(1)} $\chigraph{\R^2}{C_4}\le 4$,  \ \textup{(2)} if  $\chi(\R^2)=7$ then $\chigraphind{\R^2}{C_4}>2$, and \  \textup{(3)} $\chigraphind{\R^3}{C_4}>2$.
\end{proposition}

\subsection{Background on Euclidean Ramsey theory} \label{sec:background}

Recall that $\chi(\R^n)$ denotes the smallest $r\in \N$ such that there exists an $r$-coloring of $\R^n$ with no two points at distance $1$ and of the same color.
Despite more than $70$ years of study, the problem of determining $\chi(\R^n)$ remains unsolved even in the plane. The previously best-known bounds are $5\leq\chi(\R^2)\leq7$. The upper bound comes from a colored hexagonal grid and the lower bound is due to de Grey~\cite{degrey}, see also Exoo and Ismailescu~\cite{EI20} and the book~\cite{Soifer} by Soifer for a history on this problem. 
Recently, Sokolov and Voronov~\cite{SV25} showed that every $6$-coloring of the plane in which each color class is a disjoint union of polygons contains a monochromatic pair of points at distance $1$. Building on their work, OpenAI~\cite{AI26} has announced a proof that $\chi(\R^2)\geq6$. If correct, this would narrow the possible values of $\chi(\R^2)$ to $6$ and $7$.
In a different regime where $n$ grows, the currently best bounds on $\chi(\R^n)$ are
\begin{equation} \label{chiRn}
     (1.239...+o(1))^n \le \chi(\R^n) \le (3+o(1))^n.
\end{equation}
Here, the upper bound is due to Larman and Rogers~\cite{LR72}, see also Prosanov~\cite{Pros20} for an alternative proof, while the lower bound is due to Raigorodskii~\cite{Rai00}, who improved the first exponential lower bound $\chi(\R^n) \ge (1.207...+o(1))^n$ due to  Frankl and Wilson~\cite{FW81}.

As mentioned in \Cref{sec:1.1}, a set $A \subseteq \R^d$ is called \emph{Ramsey} if for every $r$, there exists $n \in \N$ such that $\chi(\R^n,A)>r$, where $\chi(\R^n,A)$ is defined as the smallest $r\in \N$ such that there exists an $r$-coloring of $\R^n$ with no monochromatic copy of $A$. Here the word \emph{copy} stands for an \emph{isometric copy}, i.e., a subset $A'\subseteq \R^n$ such that there is a bijection between $A$ and $A'$ that preserves the Euclidean distance. Some authors consider \emph{congruent copies} instead, namely, images of $A$ under the motions (isometries) of $\R^n$, where $\R^d$, the ambient space for $A$, is naturally treated as a subspace of $\R^n$. However, these are equivalent definition  since a local isometry between any two subsets of $\R^n$ can be extended to a global isometry (motion) of $\R^n$, see, e.g.,~\cite[Section~38, Property~IV]{Blum70} or \cite[Theorem~11.4]{WW12}.

It is known that a vertex set of a non-degenerate simplex or a regular polytope is Ramsey, see~\cite{FR90, Kriz1991, Cant2007}. Other specific Ramsey sets were found in~\cite{Beh25, LRW12, Kriz1992}. Moreover, the Cartesian product of two Ramsey sets is also Ramsey~\cite[Theorem~20]{EGMRSS1}. Frankl and R\"odl~\cite{FR90} proved that if $A$ is a vertex set of a non-degenerate simplex or a hyperbox, then $\chi(\R^n,A)$  grows exponentially with $n$.
From the other direction, we know that every Ramsey set must be finite~\cite[Theorem~19]{EGMRSS2} and spherical~\cite[Theorem~13]{EGMRSS1}, i.e., lie on the surface of a hypersphere. One of the biggest open problems in the field is Graham's \$$1000$ conjecture~\cite{G17}, which states that these two necessary conditions are also sufficient. We also refer the reader to a ``rival conjecture'' of Leader, Russell, and Walters~\cite{LRW12}, which states that a set is Ramsey if and only if it is a subset of a finite transitive set, see also the discussions in~\cite{Beh25, ILW, KK20, Kar22, LRW11}. Here we say that a set $A \subseteq \R^d$ is \emph{transitive} if it has a transitive symmetry group, i.e., if for all $a_1,a_2 \in A$ there is an isometry of $A$ that maps $a_1$ to $a_2$. Rather surprisingly, both conjectures have recently been disproved. P\'alv\"olgyi~\cite{DP26} constructed a $7$-point counterexample to Graham's conjecture~\cite{G17}, while OpenAI~\cite{AI26'} found a $4$-point counterexample to the conjecture of Leader, Russell, and Walters~\cite{LRW12}. Moreover, OpenAI~\cite{AI26'} announced a complete characterization of Ramsey sets in terms of a tensor condition over their coordinate fields.

Variants of this problem where, instead of avoiding isometric copies of just one fixed configuration, we consider an infinite family of configurations, have also been studied. The following results of this type are particularly notable. Graham~\cite{Grah80} showed that for all integers $r, n\ge 2$, every $r$-coloring of $\R^n$ contains a monochromatic set of  $n+1$  points that form a right-angled simplex of unit volume. Here a simplex is \emph{right-angled} if all its edges adjacent to some vertex are pairwise orthogonal. Answering a question of Erd\H{o}s and Graham in the negative, Kova{\v{c}}~\cite[Theorem~4]{Kov23} proved that for every $n\in \N$, there exist $r\in \N$ and an $r$-coloring of $\R^n$ such that no set of $2^n$ vertices of an $n$-dimensional box of unit volume is monochromatic; moreover, $25$ colors are enough in the plane~\cite[Theorem~3]{Kov23}.  Voronov~\cite{Vor23} proved that for any $\varepsilon>0$, every $6$-coloring of the plane contains a monochromatic pair of points at distance between $1$ and $1+\varepsilon$. Davies~\cite{davies2024odd} and  Davies, McCarty, and Pilipczuk~\cite{davies2023prime} showed that for each $r$, every $r$-coloring of the plane contains a monochromatic pair of  points such that the distance between them is an odd integer or a prime number, respectively.

Another related line of research deals with unavoidable configurations in measurable sets, see~\cite{ACMVZ, BPT, FYW, Kov23, KovPre24, LA20}. We explicitly mention the following special case of the general result by Lyall and Magyar~\cite[Theorem~2~(i)]{LA20}, which may be interesting to compare with our results. For every finite graph $H$, there exists $n_0(H)$ such that for each integer $n \ge n_0(H)$, every set $A \subseteq \R^n$ of positive upper Banach density contains an \emph{induced $\lambda$-copy} of $H$, i.e., an induced unit-copy of $H$ scaled by the factor $\lambda$, for all $\lambda$ sufficiently large in terms of $H,n$ and $A$. 
Moreover, one can take $n_0(H)=2$ if $H$ is a forest or a hypercube~\cite{KovPre24}. However, the minimum value of $n_0(H)$ is not known already in the case $H=K_3$, see~\cite[Problem~7]{Kov23}.

Returning to the planar case, as mentioned in \Cref{sec:1.1}, Erd\H{o}s, Graham, Montgomery, Rothschild, Spencer, and Straus~\cite[Section~1]{EGMRSS1} observed that the $2$-coloring of $\mathbb{R}^2$ with color-alternating infinite strips of width $\sqrt{3}/2$ contains no monochromatic copy of the vertex set of the unit equilateral triangle. They also conjectured that this is the only $2$-coloring of the plane with this property, which was disproved by Jel\'inek, Kyn\v{c}l, Stola\v{r}, and Valla~\cite{JKSV09}. For the vertex set $T$ of any other non-equilateral triangle, they conjectured that $2$ colors are not enough, i.e., that $\chi(\R^2,T)\ge 3$. Despite a considerable attention, this conjecture was verified only for a few special classes of triangles, see \cite{G17}. In particular, the case of a degenerate triangle $T = \{0,1,2\} \subseteq \R$ was resolved only recently by Currier, Moore, and Yip~\cite{CMY24}. From the other direction, Graham conjectured \cite[Conjecture~11.1.3]{G17} that $\chi(\R^2,T)\le 3$ for any triangle $T$. Improving on the work of Aichholzer and Perz~\cite{AP19}, Mundinger, Zimmer, Kiem, Spiegel, and Pokutta~\cite{MZKSP25} constructed explicit $3$-colorings verifying this conjecture for a large family of triangles. In general, however, no better upper bound than $\chi(\R^2,T) \le \chi(\R^2) \le 7$ is known. We note that asymmetric variants of similar problems were studied in~\cite{ArmTs, ConFox, CF, CW, CsTo, CMY26, FT}.

\subsection{Basic observations on the functions  $\chigraph{\R^n}{H}$ and  $\chigraphind{\R^n}{H}$ and known results}\label{sec:basic-observations}

As in classical Euclidean Ramsey theory, we can express the problem of finding  $\chigraph{\R^n}{H}$ as a hypergraph coloring problem. Namely, let for a graph $H$ a hypergraph $\cH$ be defined on the vertex set $\R^n$ and hyperedges corresponding to all sets of vertices that are unit-copies of $H$. Then $\chigraph{\R^n}{H}= \chi(\cH)$, the classical chromatic number of a hypergraph, that is the smallest number of colors used on vertices such that there is no monochromatic hyperedge.
Similarly, if the hyperedges of $\cH$ correspond to induced unit-copies of $H$ in $\R^n$, then $\chigraphind{\R^n}{H}= \chi(\cH)$. We shall be repeatedly using the hypergraph version of the well-known de Bruijn–Erd\H{o}s theorem~\cite{debruijn} on infinite graphs, which relies on the axiom of choice. Applied to our problem, this theorem can be stated as follows.

\begin{lemma}[De Bruijn--Erd\H{o}s theorem~\cite{debruijn}]
\label{debrujin}
For any finite unit-distance graph $H$ in $\R^n$, there is a finite unit-distance graph $G$ in $\R^n$ such that $G \ind{r} H$, where $r=\chigraphind{\R^n}{H}-1$. The same holds for the non-induced case.
\end{lemma}

We shall also be repeatedly using the following result on products of unit-distance graphs.

\begin{lemma}[{Horvat--Pisanski~\cite[Theorem~3.4]{HP10}}]
\label{product}
For any integer $n\geq2$, if $G$ and $H$ are unit-distance graphs in $\R^n$, then $G\square H$ is a unit-distance graph in $\R^n$. In particular, if $G$ is a unit-distance graph in $\R^n$, then $G^{\square N}$ is a unit-distance graph in $\R^n$ for each $N \in \N$.
\end{lemma}

We proceed with basic observations. First, it is immediate from the definitions that a graph $H$ has at least one edge, then $\chigraphind{\R^n}{H}\leq\chigraph{\R^n}{H}$.
If $G$ is a subgraph (resp. an induced subgraph) of $H$, then 
$\chigraph{\R^n}{H}\leq \chigraph{\R^n}{G} \  (\text{resp. }\chigraphind{\R^n}{H}\leq \chigraphind{\R^n}{G}).$
Moreover, if a set $A\subseteq\R^d$ contains an induced unit-copy of $H$ as a subset, then $\chi(\R^n, A) \le \chigraphind{\R^n}{H}$ for all $n$.

Note that cliques are rigid as unit-distance graphs in the following sense. Let $\triangle_{d}$ be a $d$-element vertex set of a $(d-1)$-dimensional unit regular simplex. Then every subset $A \subseteq \R^n$ is a unit-copy of $K_{d}$ if and only if $A$ is an isometric copy of $\triangle_{d}$. Hence $\chigraphind{\R^n}{K_d}= \chigraph{\R^n}{K_d} = \chi(\R^n,\triangle_d)$.
The best known bounds for simplices are $(\psi^{1/d}+o(1))^n \le \chi(\R^n,\triangle_d) \le (1+\sqrt{2d/(d-1)}+o(1))^n$, as $n \to \infty$,  where $\psi=1.239\dots$ is the constant from~\eqref{chiRn}. Here, the lower bound is due to Kupavskii, Zakharov, and the third author~\cite{KSZ22}, while the the upper bound is due to Prosanov~\cite[Corollary~2]{Pros18}.

Note that for every graph $H$ on $m$ vertices with clique number $\omega$, we have  $\chi(\R^n,\triangle_m) \le \chigraph{\R^n}{H} \le  \chi(\R^n,\triangle_\omega).$ In particular, we have
$\chigraph{\R^n}{H}\leq \chi(\R^n,\triangle_2) = \chi(\R^n) \leq\left(3+o(1)\right)^n.$
Besides, $\chigraph{\R^n}{H}$ grows exponentially with $n$ for every fixed graph $H$. Similar statement for the induced version is less straightforward, see \Cref{thm:main-asymptotic}~(1) below.

\subsection{Results on $\chigraph{\R^n}{H}$ and $\chigraphind{\R^n}{H}$ for growing dimensions} \label{sec:growing}

Next, we state our asymptotic results on growing dimensions. Note that all the constants produced by our proof are explicit.

\begin{theorem} \label{thm:main-asymptotic}
Let $H$ be a graph. The following holds as $n \to \infty$: 
   \begin{enumerate}
   \setlength\itemsep{1ex}
   \item There is $c=c(H)>1$ such that $\chigraphind{\R^n}{H} \ge  (c+o(1))^n$.
   \item For each $m\ge2$, there is $c_m>1$ such that if $H$ is $m$-partite, then  $\chigraph{\R^n}{H} \ge (c_m+o(1))^n$. 
   Moreover, for $m=2$, one can take $c_2= (4/3)^{1/4}\sim 1.074\dots$.
   \item For each $\ell\ge 1$, we have $\chigraph{\R^n}{C_{2\ell+1}} \le (1+ 2\cos(\frac{\pi}{4\ell+2})+o(1))^n$. In particular, if $H$ is not bipartite, then there is some $\varepsilon=\varepsilon(H)>0$ such that $\chigraph{\R^n}{H} \le (3-\varepsilon + o(1))^n$.
    \end{enumerate}
\end{theorem}

Recently, Cheng and Xu~\cite{CX} initiated a study of ``canonical'' problems in Euclidean Ramsey theory, where the number of colors is not specified, and one is  looking for ``unavoidable'' colorings of copies of a fixed configuration, usually monochromatic or rainbow ones. Recall that a coloring of a set is \emph{rainbow} if each element gets a unique color. Geh\'er, Sagdeev, and T\'oth~\cite[Theorem~4]{GST} showed for any $d\in \N$, there is a sufficiently large $n$ such that every coloring of $\R^n$ contains a copy of $\{0,1\}^d$ that is either monochromatic or rainbow. Fang, Ge, Shu, Xu, Xu, and Yang~\cite{FGSXXY} proved a similar result for 
vertex sets of rectangles, non-degenerate triangles, and 
special classes of simplices. In this paper, we use the aforementioned result from~\cite{GST} to deduce the following canonical result for graphs in Euclidean Ramsey theory.

\begin{theorem} \label{thm:canon}
    For each graph $H$, there is a sufficiently large $n$ such that every coloring of $\R^n$ contains an induced unit-copy of $H$ that is either monochromatic or rainbow.
\end{theorem}

The rest of the paper is structured as follows.
We prove our results on Ramsey properties of Cartesian powers, i.e., \Cref{thm:Turan}, \Cref{pro:Turan}, and \Cref{lem:layered} in Section \ref{sec:cartesian}. 
In Section~\ref{sec:euclidean} we prove our main result on graphs in Euclidean Ramsey theory, specifically \Cref{thm:main} and \Cref{general}.  Section \ref{small} is devoted to Ramsey properties of some small graphs in Euclidean spaces, in particular to \Cref{thm:main-layer} and \Cref{leftovers}.
In Section \ref{sec:growing-canonical} we address the growing dimensions and the canonical-type statement. 
Finally in Section \ref{sec:conclusions} we give concluding remarks and state many open problems. 
In \Cref{app:hypercube-Turan} we show an induced variant of a zero hypercube Tur\'an density theorem. In \Cref{app:exponential-lb} we present three different techniques for Euclidean Ramsey results in growing dimensions.  

\noindent
\textbf{Note added.} Since the first version of this paper appeared online, several of our results have been strengthened and extended. Alm\'asi, Axenovich, and Sagdeev~\cite{AlAxSag} proved that $\chigraph{\R^n}{H} = \chi(\R^n)$ for any layered graph $H$ and any $n\ge 2$. In particular, this extends \Cref{thm:main}~(2) to $\ell=3,5$ and strengthens \Cref{general}~(2) and \Cref{thm:main-layer}~(1). Moore and Sagdeev~\cite{MS} showed that $\chigraphind{\R^2}{C_4}>2$, which improves \Cref{leftovers}~(2). Kula, Omar, Stockwell, and Xie~\cite{KOSX} considered the case of multiple forbidden distances and quantitatively improved the value of $c(H)$ from \Cref{thm:main-asymptotic}~(1) for some graphs $H$. These results also answer some of the questions we pose in Section~\ref{sec:conclusions}, see the respective remarks there.

\section{Proofs of  Ramsey-type results on Cartesian powers} \label{sec:cartesian}

\subsection{Preliminary lemmas}

\begin{lemma}
\label{zero_density}
\
\begin{enumerate}
    \setlength\itemsep{1ex}
    \item If $H$ is a  vertex-transitive graph, then every $H$-forest has zero induced $H$-slice density.
    \item Let $H$ be a graph and $u,v$ be two vertices of $H$. Let $\Gamma$ be a bipartite graph with parts $A$ and $B$ of zero strong hypercube Tur\'an density with respect to $(A,B)$. Then the graph $\Gamma_{A,B}(H,u,v)$ has zero $H$-slice density.
\end{enumerate}
\end{lemma}

\begin{proof}[Proof of~\Cref{zero_density} (1)]
Let $F$ be an $H$-forest. We can assume by adding some copies of $H$ if needed that $F$ is an \emph{$H$-tree}, that is a union of pairwise edge-disjoint copies $H_1,\dots,H_m$ of $H$ such that every $V(H_i)$ intersects $\bigcup_{j<i}V(H_j)$ in exactly one vertex. Fix any constant $\varepsilon>0$ and let $N=N(F,\varepsilon)\in\N$ be sufficiently large.

Recall that an \emph{$H$-slice} in $H^{\square N}$ is a subgraph induced by a Cartesian product of $V(H)$ and $N-1$ one-element subsets of $V(H)$ in some order. Note that every $H$-slice is an induced copy of $H$. We call the coordinate, in which the vertices of an $H$-slice vary, the \emph{direction} of that $H$-slice. For example, the direction of the $H$-slice induced by $V(H)\times\{x\}\times\{y\}\times\dots$ is $1$. Note that each vertex of $H^{\square N}$ belongs to exactly $N$ $H$-slices of pairwise distinct directions.

Let $G$ be a subgraph of $H^{\square N}$ with at least $\varepsilon$-fraction of $H$-slices. On average, each vertex of $G$ is contained in at least $\varepsilon N$ $H$-slices that are subgraphs of $G$.
By iteratively removing vertices that are contained in strictly less than $\varepsilon N/|V(H)|$ such $H$-slices, we obtain a non-empty induced subgraph $G'$ of $G$ such that every vertex of $G'$ is contained in at least $\varepsilon N/|V(H)|$ $H$-slices that are subgraphs of $G'$.

We shall find an induced copy of $F$ in $G'$ by greedily adding one $H$-slice with a new direction after another. Formally, for each $t\in[m]$, we find an induced subgraphs $F_t$ in $G'$, where $F_t$ is a copy of $\bigcup_{i=1}^tH_i$, such that the vertices of $F_t$ coincide in all but at most $t$ coordinates. Then $F_m$ would be the desired induced copy of $F$.

Let $F_1$ be an arbitrary $H$-slice from $G'$, let $i_1$ be the direction of $F_1$ and $I_1=\{i_1\}$. Assume that $F_t$ has been embedded as desired, $t<m$.
Let $v$ be the vertex in $F_t$ that corresponds to the unique vertex shared by $\bigcup_{i=1}^tV(H_i)$ and $V(H_{t+1})$. Let $I_t\subseteq[m]$ denote the set of coordinates in which the vertices of $F_t$ do not agree. Consider an $H$-slice $H'$ in $G'$ that contains $v$ and has direction $i_{t+1}$ not in $I_t$. Such $H'$ exists since the number of $H$-slices in $G'$ containing $v$ is at least $\varepsilon N/ |V(H)| > m > t$.
Let $F_{t+1}= F_t \cup H'$ and $I_{t+1} = I_t\cup \{i_{t+1}\}$. Since $H$ is vertex-transitive, $H'$ is a copy of $H$, and $V(H')\cap V(F_t)=\{v\}$, we have that $F_{t+1}$ is a copy of $\bigcup_{i=1}^{t+1}H_i$.

We need to argue that $F_{t+1}$ is an induced subgraph of $G'$.
Since $F_t$ and $H'$ are induced subgraphs of $G'$,  it is sufficient to argue that there are no edges ``between" $F_t$ and $H'$ except those incident to $v$.
Consider a  vertex $w'\in V(H')\sm \{v\}$ and any vertex $w\in V(F_t)\sm \{v\}$.
We see that $w'$ differs from $v$ in the $i_{t+1}$st coordinate and thus $w'$ differs from $w$ in the $i_{t+1}$st coordinate. In addition, $w'$ coincides with $v$ in all other coordinates and $v$ differs from $w$ in some other, $j$th coordinate, $j\in I_t$. Thus $w'$ and $w$ differ in the $j$th coordinate. This implies that $w'$ and $w$ differ in at least two coordinates and thus are not adjacent in $H^{\square N}$.
\end{proof}

\begin{proof}[Proof of~\Cref{zero_density} (2)]
Fix any $\varepsilon>0$. Let $d=d(\Gamma,\varepsilon)\in\N$ and $N=N(H,d,\varepsilon)\in\N$ be sufficiently large. Let $G$ be a subgraph of $H^{\square N}$ containing at least $2\varepsilon N |V(H)|^{N-1}$ $H$-slices. We shall show that if $\Gamma$ has zero strong hypercube Tur\'an density with respect to $(A,B)$, then $G$ contains a copy of $\Gamma_{A,B}(H,u,v)$.

Let $F$ be the graph on $V(H^{\square N})$, in which two vertices are adjacent if and only if one is obtained from the other by replacing one of the coordinates equal to $u$ with $v$. Observe that each connected component of $F$ is a copy of $Q_t$ for some $t\leq N$ as its vertex set is a Cartesian product of $t$ copies of the set $\{u,v\}$ and $N-t$ one-element subsets of $V(H)\sm\{u,v\}$ in some order. In particular, there are precisely $\binom{N}{t}(|V(H)|-2)^{N-t}$ components that are copies of $Q_t$, and for each $t$ with $1\leq t<d$, the total number of edges in them is at most
\begin{equation*}
    |E(Q_t)|\binom{N}{t}(|V(H)|-2)^{N-t}\leq N^{d}(|V(H)|-2)^{N-1}<\frac{\varepsilon}{d}N|V(H)|^{N-1},
\end{equation*}
provided that $N$ is sufficiently large. In the last inequality above we used the fact that $H$ is fixed, namely, $|V(H)|-2=(1-\delta)|V(H)|$ for some constant $\delta>0$. Since $|E(F)|=N|V(H)|^{N-1}$, the components that are copies of $Q_t$ for some $t<d$ contain in total at most $\varepsilon|E(F)|$ edges.

Let $F'$ be the spanning subgraph of $F$ consisting of all edges whose endpoints belong to an $H$-slice in $G$. Since $G$ contains at least $2\varepsilon N |V(H)|^{N-1}$ $H$-slices and each $H$-slice contributes exactly one edge to $F'$, we have $|E(F')|\geq2\varepsilon N |V(H)|^{N-1}=2\varepsilon|E(F)|$. Thus $F'$ contains at least an $\varepsilon$-fraction of the edges of some connected component $Q$ of $F$ that is isomorphic to $Q_t$ for some $t\geq d$.

Recall that $V(Q)$ is a Cartesian product of $t$ copies of the set $\{u,v\}$ and $N-t$ one-element subsets of $V(H)\sm\{u,v\}$. We identify $V(Q)$ with $V(Q_t)$ by ignoring the $N-t$ fixed coordinates and identifying $u$ with $1$ and $v$ with $0$.

Since $d$ is sufficiently large, $\Gamma$ has zero strong hypercube Tur\'an density with respect to $(A,B)$, and $F'$ has at least an $\varepsilon$-fraction of the edges of $Q$, we conclude that $F'$ contains a copy $\Gamma'$ of $\Gamma$. Moreover, $\Gamma'$ lies within one edge layer of $Q$, with $A'$ above $B'$, where $A'$ and $B'$ are the parts of $\Gamma'$ corresponding to $A$ and $B$, respectively.

Note that every edge $xy\in E(\Gamma')$, $x\in A'$ and $y\in B'$, is an edge in $F'$, where $y$ is obtained from $x$ by switching one coordinate from $u$ to $v$. Namely, every edge $xy\in E(\Gamma')$ corresponds to a unique $H$-slice in $G$, where $x$ and $y$ correspond to the vertices $u$ and $v$, respectively. Furthermore, for any two non-adjacent edges $xy,x'y'\in E(\Gamma')$, we claim their corresponding $H$-slices are vertex-disjoint. Indeed, recall that the direction of an $H$-slice is the index of the coordinate where its vertices differ. Let $H_1$ and $H_2$ be the corresponding $H$-slices of $xy$ and $x'y'$, respectively. Suppose $H_1$ and $H_2$ share some vertex $w$. Then $H_1$ and $H_2$ must have different directions $i$ and $j$, as otherwise we have $H_1=H_2$ and thus $xy=x'y'$, a contradiction. Accordingly, the $i$th coordinate of $w$ is equal to the $i$th coordinate of $x'$ which is in $\{u,v\}$. This particularly implies that $w$ has the role of either $u$ or $v$ in $H_1$, yielding that $\{x,y\}\cap \{x',y'\}\neq\varnothing$, a contradiction. Therefore, the subgraph $\Gamma'$ in $F'$ corresponds to a subgraph $\Gamma'_{A',B'}(H,u,v)$ in $G$, which is a copy of $\Gamma_{A,B}(H,u,v)$.
\end{proof}

Let $K_3$ be a complete graph on the vertex set $[3]$. Note that two vertices from $[3]^N$ are adjacent in $K_3^{\square N}$ if any only if they differ in exactly one coordinate.
For each $\x=(x_1,\dots,x_N)\in [3]^N$, we let $I(\x)$ be the set of indices $i$ such that $x_i\neq 3$. For $S \subseteq [N]$, let $\x|_S$ be a tuple obtained from $\x$ by deleting all coordinates with indices not in $S$. For example, $I \coloneq I((2, 1, 3, 1, 3, 1))=\{1,2,4,6\}$ and $(2, 1, 3, 1, 3, 1)|_I = (2, 1, 1, 1)$.
We say that two elements $\x$ and $\y$ from $[3]^N$ are \emph{equivalent}, and write $\x\sim \y$, if $\x|_{I(\x)} = \y|_{I(\y)}$. For $I=I(\x)$, we have  $\x|_{I} \in [2]^{|I|}$. For example, $(2,1,3,1,3,1) \sim (3, 2, 1, 3, 1, 1)$ and $(2,1,3,1,3,1) \not\sim (3, 1, 3, 1, 1, 2)$. Note that if two distinct elements $\x$ and $\y$ are equivalent then $|I(\x)| =|I(\y)| \neq 0$.

\begin{lemma}
\label{lem:equivalent}
For any $d' \in \N$, there exists $N = N(d')$ such that the following holds. For every $2$-coloring of $[3]^{N}$, there is a $d'$-element subset $S' \subseteq[N]$ such that every two equivalent elements $\x, \y \in [3]^N$ with $I(\x)\subseteq S',\, I(\y)\subseteq S'$ receive the same color.
\end{lemma}

\begin{proof}[Proof of~\Cref{lem:equivalent}]
The proof is a relatively straightforward iteration of the Ramsey theorem. The analogous statement for $[2]^N$ instead of $[3]^N$ is a known Layered Lemma, see for example \cite{AW, AL}. A similar argument was also used by Leader, Russell, and Walters as part of their proof of~\cite[Theorem~3.1]{LRW12} and by Ivan, Leader, and Walters~\cite[Theorem~3]{ILW}. Nevertheless, we give a full proof here for completeness.

Let $N=N(d')$ be sufficiently large.
For a tuple $\aaa \in [2]^{k}$, $1\leq k<N$, and $I\subseteq [N]$ of size $k$, we let $\x^I[\aaa]=\x$ denote the unique tuple in $[3]^{N}$ such that $I(\x)=I$ and $\x|_I = \aaa$. Namely, $\x^I[\aaa]$ is obtained by filling the designated positions $I$ with the tuple $\aaa$ and filling the remaining positions with $3$. For example, if $N=5$, then $\x^{\{2,4\}}((1,1))=(3, 1, 3, 1, 3)$.

Consider an arbitrary $2$-coloring $c$ of $[3]^{N}$ with colors $r$ and $b$ standing for red and blue. To find the desired subset $S'\subseteq[N]$, we iteratively apply the Ramsey theorem to certain auxiliary hypergraphs $\cH_1, \dots, \cH_{d'}$ with a new auxiliary edge-coloring defined on each of them. Recall that the Ramsey theorem claims that for any $r, q, p$ there is a sufficiently large $n$ such that any $r$-coloring of a complete $q$-uniform hypergraph on $n$ vertices results in a complete monochromatic subhypergraph on $p$ vertices. 

{\bf Step 1.} Let $\cH_1$ be the $1$-uniform complete hypergraph on the vertex set $[N]$. We assign a color given by the set $\{(a,c(\x^{\{i\}}[a])):a\in[2]\}$ to each edge $\{i\}$ of $\cH_1$. For example if  $c((3, 1, 3, \ldots, 3))=r$ and  $c((3, 2, 3, \ldots, 3))=b$, then the color of $\{2\}$ is $\{(1, r), (2,b)\}$. This gives an auxiliary coloring of the edges of $\cH_1$ using at most $4$ colors. By the pigeonhole principle, there is a large set of vertices of $\cH_1$, say $S_1\subseteq [N]$, inducing a monochromatic $1$-uniform clique under the auxiliary coloring. Observe that any two equivalent elements $\x, \y \in [3]^N$ such that $I(\x)\subseteq S_1,\, |I(\x)|=1$ and $I(\y)\subseteq S_1, \, |I(\y)|=1$  satisfy $c(\x)=c(\y)$. 

Let $k$ be an integer with $2\leq k\leq d'$. Assume steps $1, \dots, k-1$ have been performed and a subset $S_{k-1} \subseteq[n]$ has been found with the property that any two equivalent elements $\x, \y \in [3]^N$ such that $I(\x)\subseteq S_{k-1},\, |I(\x)|\le k-1$ and $I(\y)\subseteq S_{k-1}, \, |I(\y)|\le k-1$  satisfy $c(\x)=c(\y)$. 
 
{\bf Step k.} Let $\cH_k$ be the $k$-uniform complete hypergraph on the vertex set $S_{k-1}$. We assign to each edge $E\in\binom{S_{k-1}}{k}$ of $\cH_k$ a color given by the set $\{(\aaa,c(\x^{E}[\aaa])):\aaa\in[2]^k\}$. This auxiliary edge-coloring of $\cH_k$ uses at most $2^{2^k}\leq2^{2^{d'}}$ colors. Since $S_{k-1}$ is large enough (as $N=N(d')$ was initially chosen to be sufficiently large), by the hypergraph Ramsey theorem, there is a large subset $S_k$ of $S_{k-1}$ inducing a monochromatic $k$-uniform clique in $\cH_k$ under the auxiliary edge-coloring. 
Observe that any two equivalent elements $\x, \y \in [3]^N$ such that $I(\x)\subseteq S_k,\, |I(\x)|=k$ and $I(\y)\subseteq S_k, \, |I(\y)|=k$  satisfy $c(\x)=c(\y)$.  In addition, since $S_k \subseteq S_{k-1}$, we also have $c(\x)=c(\y)$ for any two equivalent elements $\x, \y \in [3]^N$ such that $I(\x)\subseteq S_k,\, |I(\x)|\le k-1$ and $I(\y)\subseteq S_k, \, |I(\y)|\le k-1$.

After $d'$ steps, we obtain a set $S_{d'}\subseteq [N]$ with the property that $|S_{d'}|\geq d'$ and any two equivalent elements $\x, \y \in [3]^N$ such that $I(\x)\subseteq S_{d'},\, |I(\x)|\le d'$ and $I(\y)\subseteq S_{d'}, \, |I(\y)|\le d'$  satisfy $c(\x)=c(\y)$. The desired set $S'$ is obtained by taking any $d'$-element subset of $S_{d'}$.
\end{proof}

\subsection{Proofs of Theorem 1.4, Theorem 1.2, and Proposition 1.3}

\begin{proof}[Proof of~\Cref{lem:layered}]
 Let $H$ be an (induced) subgraph of the $k$th edge layer of $Q_d$. For $d'=dk+d-1$, let $N=N(d')$ be from the statement of \Cref{lem:equivalent}. Consider an arbitrary $2$-coloring $c$ of $[3]^{N}$. By \Cref{lem:equivalent}, there is a $d'$-element subset $S' \subseteq[N]$ such that every two equivalent elements in $Q \coloneq \{\x \in [3]^N: I(\x)\subseteq S'\}$ receive the same color.

 First, suppose that there are elements $\x, \z \in Q$ of the same color such that exactly $d-k$ coordinates of $\x|_{S'}$ are equal to $3$, $\x|_{I(\x)}$ contains a block $B$ of $k$ consecutive $1$'s, and $\z$ is obtained from $\x$ by replacing one of these $1$'s with 3. Since equivalent elements in $Q$ receive the same color, by picking a suitable element from the equivalence class, we can assume without loss of generality that $\x|_{S'}$ is obtained from $\x|_{I(\x)}$ by inserting $3$'s in the block $B$ or next to it. Let $S$ be the set of indices $i\in S'$ such that either $x_i=3$ or $x_i=1$ for some $1$ from $B$. Note that $|S|=k+(d-k)=d$.
 For example, if $k=2$, $d=3$, $N=9$, $S'=[9]\sm\{6\}$, $\x=(1,2,2,2,1,3,1,3,2)$, and $\z=(1,2,2,2,3,3,1,3,2)$, then $\x|_{I(\x)} = (1,2,2,2,\textbf{1},\textbf{1},2)$, where the block $B$ is in bold,  $\x|_{S'} = (1,2,2,2,1,1,3,2)$ is obtained from $\x|_{I(\x)}$ by inserting 3 to the right of $B$, and   $S=\{5,7,8\}$.
 
 Consider the set $L$ of elements $\y \in [3]^N$ such that $\y|_{[N]\sm S} = \x|_{[N]\sm S}$ and either $k$ or $k-1$ coordinates of $\y|_S$ are equal to $1$, while the remaining ones are equal to $3$. Observe that each $\y \in L$ is equivalent to either $\x$ or $\z$ and thus the set $L$ is monochromatic. Moreover, the subgraph of $K_3^{\square N}$ induced by $L$ is a copy of the $k$th edge layer of $Q_d$ and thus contains an (induced) copy of $H$, as desired. In the example above, $L$ consists of 6 elements of the form $(1,2,2,2,*,3,*,*,2)$ where either one or two $*$'s are equal to $1$, and the remaining $*$'s are equal to $3$. Each $\y \in L$ is equivalent to either $\x$ (if two of its $*$'s are equal to $1$) or to $\z$ (otherwise). Moreover, by deleting the fixed coordinates, i.e., the coordinates whose indices are not in $S$, and identifying each $*$ equal to $3$ with $0$, we see that the subgraph of $K_3^{\square 9}$ induced by $L$ is a copy of the $2$nd edge layer of $Q_3$, as claimed.

 Second, suppose that there are elements $\y, \z \in Q$ of the same color such that exactly $d-k$ coordinates of $\y|_{S'}$ are equal to $3$, $\y|_{I(\y)}$ contains a block of $k$ consecutive $2$'s and $\z$ is obtained from $\y$ by replacing one of these $2$'s with $3$. In this case, a similar argument leads to an (induced) monochromatic copy of $H$, as desired.

 Now consider the set $Q'$ of all $2^d$ elements $\x \in Q$ such that, depending on the parity of $d$, $\x|_{S'}$ is of the form $(\textbf{1}_{k-1},*,\textbf{2}_{k-1},*, \cdots ,*,\textbf{1}_{k-1},\textbf{3}_{d-k})$ or $(\textbf{1}_{k-1},*,\textbf{2}_{k-1},*,\cdots ,*,\textbf{1}_{k-1},*,\textbf{2}_{k-1},\textbf{3}_{d-k})$ where $\textbf{1}_{k-1}$ denotes a block of $1$'s of length $k-1$, $\textbf{2}_{k-1}$ denotes a block of $2$'s of length $k-1$, these blocks alternate and there are $d+1$ of them in total, each $*$ between these blocks is either $1$ or $2$, and all $d-k$ coordinates equal to $3$ are grouped in the end.

 For the example above when $k=2$, $d=3$, $N=9$, and $S'=[9]\sm\{6\}$, the set $Q'$ consists of 8 elements of the form $(1,*,2,*,1,3,*,2,3)$, where each $*$ is either $1$ or $2$.

 Consider any two elements $\x,\y \in Q'$ that differ in exactly one coordinate, say $x_{i}=1, y_{i}=2$. Let $\z$ be the element of $Q$ such that $z_{i}=3$ while all the other coordinates of $\z$ coincide with those of $\x$ and $\y$. Note that the vertices $\x, \y, \z$ are pairwise adjacent in $K_3^{\square N}$, and thus some two of them are of the same color by the pigeonhole principle. Since $\x|_{I(\x)}$ contains a block of $k$ consecutive $1$'s, the case $c(\x)=c(\z)$ leads to the desired (induced) monochromatic copy of $H$ as we have already seen. Similarly, since $\y|_{I(\y)}$ contains a block of $k$ consecutive $2$'s, the case $c(\y)=c(\z)$ leads to the same conclusion. Hence, the only case left to consider is when $c(\x)=c(\y)$ for every two elements $\x,\y \in Q'$ that differ in exactly one coordinate, i.e., when $Q'$ is monochromatic. In this case, by deleting the fixed coordinates, i.e., the coordinates not labeled by $*$'s in the definition of $Q'$, we see that the subgraph of $K_3^{\square N}$ induced by $Q'$ is a copy of $Q_d$ and thus contains the desired (induced) monochromatic copy of $H$ as well. This concludes the proof of~\Cref{lem:layered}.
\end{proof}

\begin{proof}[Proof of~\Cref{thm:Turan}]
Let $G_1,\dots,G_m$, $H_1, \ldots, H_m$, and $G$ be graphs and $r$ be a positive integer, let $V=V(G)$. 
Let $H_i$ have zero (induced) $G_i$-slice density for all $i \in [m]$. Suppose that every $r$-coloring of the vertices of $G$ contains a monochromatic (induced) copy of $G_i$ for some $i \in [m]$.
For $i \in [m]$, consider all (induced) copies of $G_i$ in $G$ and label the vertex sets of these copies by $V_i(1), \ldots, V_i(c_i)$.  We say that a copy of $G_i$ in $G$ has index $j$ if it is induced by $V_i(j)\subseteq V$.
Let $c \coloneq c_1+\dots+c_m$.  Fix $\varepsilon = \frac{1}{2cr}$ and $d_0\in \N$ such that for every integer $d \ge d_0$, any subgraph of $G_i^{\square d}$ that contains at least an $\varepsilon$-fraction of its slices contains an (induced) copy of $H_i$, $i=1, \ldots, m$.
Let $N = N(d_0, |V|)\in \N$ be sufficiently large. 

Consider an $r$-coloring of the vertices of $G^{\square N}$, i.e., of $V^N$. Recall that each of the $N|V|^{N-1}$ $G$-slices of $G^{\square N}$ contains a monochromatic (induced) copy of $G_i$ for some $i \in [m]$. 
Fix one such copy in every $G$-slice.  Each such copy has some color,  is isomorphic to some $G_i$, and has some index $j$, $1\leq j\leq c_i$. 
By the pigeonhole principle, there is a tuple $(col, i, j)$, say without loss of generality $(col, i,j)=(red, 1, 1)$, such that in $G^{\square N}$, there are at least $\frac{1}{cr}N|V|^{N-1}$ red copies of $G_1$ of index $1$. Let $V'= V_1(1)$.

Let $\cF$ be the family of all Cartesian products of $V'$ and $N-1$ one-element subsets of $V$ taken in all possible orderings of their factors, and let $\cF_{red}$ be the family of red members of $\cF$.
Since  $|\cF|= N|V|^{N-1}$ and  $|\cF_{red}|\geq \frac{1}{cr}N|V|^{N-1}$, we have that   $|\cF_{red}|\geq \frac{1}{cr}|\cF|=2\varepsilon|\cF|$. From now on, we shall only consider the vertex sets of Cartesian products and will come back to the graphs induced by them at the end of the proof.

We shall decompose $V^N$ into families of special Cartesian products.
For each $d=0, \ldots, N$, let $\cV_d$ be the family of all Cartesian products of $d$ copies of $V'$ and $N-d$ one-element subsets of $V\sm V'$, taken in all possible orderings of their factors. 

The members of all $\cV_d$, $d=0, \ldots, N$, are pairwise disjoint and form a partition of $V^N$.
Moreover, any member of $\cF$ is a subset of some member of $\cV_d$, for some $d=1, \ldots, N$. We call a member of $\cF$ \emph{bad} if it is a subset of some member of $\cV_d$, for some $d<d_0$.

First, we shall argue that the family $\cF_{bad}$ of all bad members of $\cF$ is small. Indeed, since $|V\sm V'| < |V|$, we have
\begin{equation*}
    |\cF_{bad}|=\sum_{d=1}^{d_0-1} d|V'|^{d-1}\binom{N}{d}|V\sm V'|^{N-d}<  d_0^2 |V'|^{d_0} N^{d_0}|V\sm V'|^{N}< \varepsilon N|V|^{N-1} = \varepsilon |\cF|,
\end{equation*}
provided that $N$ is sufficiently large and  treating $|V|$ and $d_0$ as fixed constants. In the last inequality we also used the fact that $|V\sm V'|\leq(1-\delta)|V|$ holds for some constant $\delta>0$.
Thus there are at least $|\cF_{red}| - |\cF_{bad}| \geq  \varepsilon|\cF|$
members of $\cF$ that are red and not bad. In particular, there is some member $Q$ of $\cV_d$, for some fixed $d\geq d_0$, such that at least an $\varepsilon$-fraction of its subsets that are members of $\cF$ are red.
Let, without loss of generality, $Q=(V')^{d}\times\{v\}^{N-d}$ for $v \in V \sm V'$.

Now, we come back to the graphs and consider the induced subgraph $F$ of $G^{\square N}$ on $Q$. We see that $F$ is isomorphic to $G_1^{\square d}$. Moreover, the vertex sets of the $G_1$-slices of $F$ are precisely those subsets of $Q$ that are members of $\cF$. Hence, at least an $\varepsilon$-fraction of $G_1$-slices of $F$ are red. Thus the subgraph of $F$ induced by the red vertices is an induced subgraph of $G^{\square N}$ that contains an (induced) red copy of $H_1$, as desired.
\end{proof}

\begin{proof}[Proof of~\Cref{pro:Turan}]
Note that $H$ has zero (induced) hypercube Tur\'an density if and only if $H$ has zero (induced) $K_2$-slice density. Moreover, from $\chi(G)=r+1$ we have that $G \xrightarrow{r} K_2$. Now, item (2) follows directly from~\Cref{thm:Turan} by letting $m=1$, $G_1=K_2$, and $H_1=H$.

Next we prove item (1). Let $H$ be a  vertex-transitive graph and $F$ be an $H$-forest. \Cref{zero_density}~(1) implies that $F$ has zero induced $H$-slice density. Then for any graph $G$ with $G\xrightarrow{r}H$ (resp., $G\ind{r}H$), by~\Cref{thm:Turan} there exists some $N\in\N$ such that $G^{\square N}\xrightarrow{r}F$ (resp. $G^{\square N}\ind{r}F$).

The proof of item (3) is essentially the same as that of item (1), once we note that if $\Gamma$ has zero strong hypercube Tur\'an density with respect to $(A,B)$, then $\Gamma_{A,B}(H,u,v)$ has zero $H$-slice density by~\Cref{zero_density}~(2).
\end{proof}

\section{Proofs of general  results on graphs in Euclidean Ramsey theory}\label{sec:euclidean}

In this section we shall prove ~\Cref{general} and~\ref{thm:main}.

\begin{proof}[Proof of~\Cref{general}]
We first prove item (1). Suppose $H$ is a  vertex transitive graph,  $F$ is an $H$-forest, and $r=\chigraphind{\R^n}{H}-1$. It is immediate that $\chigraph{\R^n}{F}\leq\chigraph{\R^n}{H}$ and $\chigraphind{\R^n}{F}\leq\chigraphind{\R^n}{H}$, since $F$ contains an induced copy of $H$. It suffices to show $\chigraph{\R^n}{F}>\chigraph{\R^n}{H}-1$ and $\chigraphind{\R^n}{F}>\chigraphind{\R^n}{H}-1$. Here we only show $\chigraphind{\R^n}{F}>\chigraphind{\R^n}{H}-1=r$, as the proof for the non-induced case is similar. By~\Cref{debrujin} we know there is a finite unit-distance graph $G$ in $\R^n$ such that $G\ind{r}H$. By~\Cref{pro:Turan}~(1) we know that $G^{\square N}\ind{r}F$ for some large $N\in \N$. Since $G^{\square N}$ is a unit-distance graph by~\Cref{product}, we have that $\chigraphind{\R^n}{F}>r$.

Next we prove item (2). We shall again only prove the induced case, as the proof of the non-induced case is analogous. Suppose $H$ is a graph of zero induced hypergraph Tur\'an density. It is clear that $\chigraphind{\R^n}{H}\leq\chigraph{\R^n}{}$. To show that $\chigraphind{\R^n}{H}>\chigraph{\R^n}{}-1=r$, we let $G$ be a finite unit-distance graph in $\R^n$ with $\chi(G)=r+1$. The existence of such $G$ is guaranteed by~\Cref{debrujin}. Then by~\Cref{pro:Turan}~(2) there is some $N\in\N$ such that $G^{\square N}\ind{r}H$. As $G^{\square N}$ is a unit-distance graph by~\Cref{product}, we conclude that $\chigraphind{\R^n}{H}>r$.

As the proof of item (3) mirrors that of item (2), we omit it here.
\end{proof}

\begin{proof}[Proof of~\Cref{thm:main}]
Since all even cycles of length $8$ or at least $12$ have zero induced hypercube Tur\'an density (see~\Cref{lem:a2} for a proof of this fact), items (1) and (2) of~\Cref{thm:main} are simply special cases of items (1) and (2) of~\Cref{general}, respectively. It suffices to prove item (3).

Let $G$ be a unit-distance graph in $\R^n$ with $\chi(G)=\chi(\R^n)$ and the smallest number of vertices. The existence of such $G$ follows from~\Cref{debrujin}, and $|V(G)|$ only depends on $n$. Let $\ell_1=\lfloor{|V(G)|/2}\rfloor$, and let $\ell\geq \ell_1+5=\ell_0$ be any integer. 

To see that $\chigraph{\R^n}{C_{2\ell+1}} \leq  \lceil \chi(\R^n)/ 2 \rceil$, we first properly color $\R^n$ with $\chi(\R^n)$ colors. Then group these colors into $\lceil \chi(\R^n)/ 2 \rceil$ classes, each containing at most two colors. The unit-distance graph induced by each such class is bipartite, and hence contains no odd cycle as a subgraph. This yields an odd-cycle-free coloring of $\R^n$ using $\lceil \chi(\R^n)/ 2 \rceil$ colors.

Let $r=\lceil \chi(\R^n)/ 2 \rceil-1$. To show $\chigraph{\R^n}{C_{2\ell+1}}>r$, we shall use the full strength of~\Cref{thm:Turan}. Observe that $G\xrightarrow{r}\{C_3,C_5,\dots,C_{2\ell_1+1}\}$, i.e., every $r$-coloring of $V(G)$ contains a monochromatic copy of $C_{2i+1}$ for some $i\in[\ell_1]$. Indeed, since $\chi(G)=\chi(\R^n)$, in every $(\lceil \chi(\R^n)/ 2 \rceil-1)$-coloring there exists a color class that induces a non-bipartite graph, which implies that there is a monochromatic odd cycle in $G$ whose length is at most $|V(G)|\leq2\ell_1+1$.

Fix any $i\in[\ell_1]$. Note that $2(\ell-i+1)\geq 12$. Let $H$ be a copy of $C_{2i+1}$ with two adjacent vertices $u,v\in V(H)$. Let $\Gamma$ be a copy of $C_{2(\ell-i+1)}$ with parts $A$ and $B$. Since $C_{2(\ell-i+1)}$ has zero hypercube Tur\'an desity (see~\cite{C} or~\Cref{lem:a2}) and $(A,B)$ is unique up to a relabeling due to the symmetry of $\Gamma$, we see that $\Gamma$ has zero strong hypercube Tur\'an density with respect to $(A,B)$, which together with~\Cref{zero_density}~(2) implies that $\Gamma_{A,B}(H,u,v)$ has zero $H$-slice density. Since $C_{2\ell+1}$ is a subgraph of $\Gamma_{A,B}(H,u,v)$ and $H$ is a copy of $C_{2i+1}$, we have that $C_{2\ell+1}$ has zero $C_{2i+1}$-slice density for each $i\in[\ell_1]$.
Now~\Cref{thm:Turan} applied with $m=\ell_1$, $G_i=C_{2i+1}$, $H_i=C_{2\ell+1}$, $i \in [\ell_1]$ yields that there is some $N\in\N$ such that $G^{\square N}\xrightarrow{r}C_{2\ell+1}$. Since $G^{\square N}$ is a unit-distance graph in $\R^n$ by~\Cref{product}, we have $\chigraph{\R^n}{C_{2\ell+1}}>r$, completing the proof.
\end{proof}

\section{Proofs of results on small graphs in Euclidean Ramsey theory}\label{small}

First, note that \Cref{thm:main-layer}~(1) is immediate from \Cref{lem:layered} since $K_3^{\square N}$ is a unit-distance graph in the plane for every $N \in \N$ by \Cref{product}.

\begin{proof}[Proof of~\Cref{thm:main-layer}~(2)]
Next we prove that $\chigraph{\R^2}{Q_{11}} =2$. It suffices to give a $2$-coloring of $\R^2$ avoiding any monochromatic unit-copy of $Q_{11}$. We shall use the $2$-coloring illustrated in~\Cref{stair_coloring}, which corresponds to alternating staircases with length $\lambda=1+3/\sqrt{2}$ and width $1$. Note that the lower boundary of each staircase is colored with the same color as the staircase. We remark that a similar coloring was used in a different context in~\cite{KS23}. 

\usetikzlibrary{patterns}
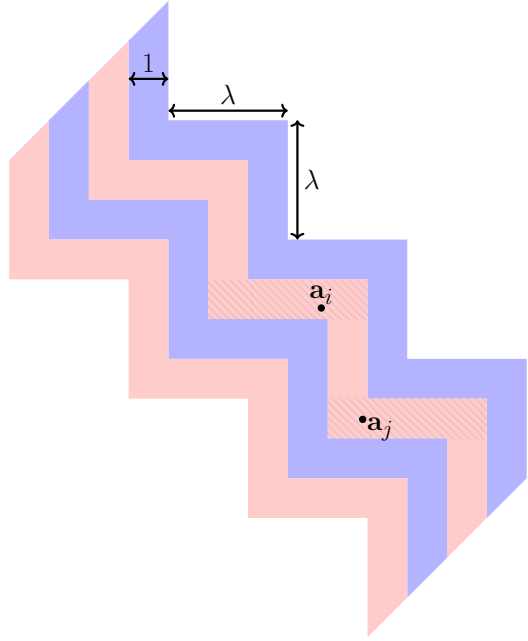
\begin{figure}[h!]
\centering
\resizebox{0.42\textwidth}{!}{%
\begin{circuitikz}

\fill[blue!30]
(2.5,20.75) -- (2.5,17) -- (6.25,17) -- (6.25,13.25) -- (10,13.25) -- (10,9.5) -- (13.75,9.5) -- (13.75,5.75) -- (12.5,4.5) -- (12.5,8.25) -- (8.75,8.25) -- (8.75,12) -- (5,12) -- (5,15.75) -- (1.25,15.75) -- (1.25,19.5) -- cycle;

\fill[red!20]
(1.25,19.5) -- (1.25,15.75) -- (5,15.75) -- (5,12) -- (8.75,12) -- (8.75,8.25) -- (12.5,8.25) -- (12.5,4.5) -- (11.25,3.25) -- (11.25,7) -- (7.5,7) -- (7.5,10.75) -- (3.75,10.75) -- (3.75,14.5) -- (0,14.5) -- (0,18.25) -- cycle;

\fill[blue!30]
(0,18.25) -- (0,14.5) -- (3.75,14.5) -- (3.75,10.75) -- (7.5,10.75) -- (7.5,7) -- (11.25,7) -- (11.25,3.25) -- (10,2) -- (10,5.75) -- (6.25,5.75) -- (6.25,9.5) -- (2.5,9.5) -- (2.5,13.25) -- (-1.25,13.25) -- (-1.25,17) -- cycle;

\fill[red!20]
(-1.25,17) -- (-1.25,13.25) -- (2.5,13.25) -- (2.5,9.5) -- (6.25,9.5) -- (6.25,5.75) -- (10,5.75) -- (10,2) -- (8.75,0.75) -- (8.75,4.5) -- (5,4.5) -- (5,8.25) -- (1.25,8.25) -- (1.25,12) -- (-2.5,12) -- (-2.5,15.75) -- cycle;

\fill[pattern=north west lines, opacity=0.3]
(3.75,12) -- (8.75,12) -- (8.75,10.75) -- (3.75,10.75) -- cycle;

\fill[pattern=north west lines, opacity=0.3]
(7.5,8.25) -- (12.5,8.25) -- (12.5,7) -- (7.5,7) -- cycle;

\draw[line width=0.7mm, <->] (2.5,17.3) -- (6.25,17.3);
\node at (4.37,17.8) {\Huge$\lambda$};

\draw[line width=0.7mm, <->] (6.55,13.25) -- (6.55,17);
\node at (7,15.13) {\Huge$\lambda$};

\draw[line width=0.7mm, <->] (1.25,18.3) -- (2.5,18.3);
\node at (1.875,18.8) {\Huge$1$};

\node [draw, shape = circle, fill = black, inner sep=2pt] at (7.3,11.1){};
\node [above] at (7.3,11.1) {\Huge$\aaa_i$};

\node [draw, shape = circle, fill = black, inner sep=2pt] at (8.6,7.6){};
\node [right] at (8.6,7.3) {\Huge$\aaa_j$};

\end{circuitikz}
}%
\caption{The staircase coloring of $\R^2$. Two red horizontal stairs are highlighted.}
\label{stair_coloring}
\end{figure}

Suppose towards a contradiction that there exists a monochromatic, say red, unit-copy $A \subseteq \R^2$ of $Q_{11}$. 
Let $\aaa_0\in A$ be a point whose orthogonal projection onto the line $x=y$ has the smallest $x$-coordinate. Due to the symmetry of $Q_{11}$, we can assume that $\aaa_0 \in \R^2$ corresponds to the vertex $\mathbf{0}\coloneqq(0, \dots, 0)$ of $Q_{11}$, and let vectors $\uu_1, \dots, \uu_{11}$ of unit length be 
such that $\aaa_0+\uu_i\in A$ correspond to the neighbors of $\mathbf{0}$ in $Q_{11}$.
Note that a $4$-element set $B\subseteq\R^2$ is a unit-copy of $Q_2$, i.e., a $4$-cycle, if and only if $B$ is the vertex set of a rhombus of side length 1. This simple observation implies that $A = \{\aaa_0+\sum_{i \in S} \uu_i: S \subseteq [11]\}$.

By the choice of $\aaa_0$, we know that each of the unit vectors $\uu_1,\dots,\uu_d$ makes an angle of at most $\pi/2$ with the vector $(1,1)$. By the pigeonhole principle and due to the symmetry of the coloring with respect to the line $x=y$, we may assume without loss of generality that six of the vectors, say $\uu_1, \dots, \uu_6$, make an angle of at most $\pi/4$ with the vector $(1,0)$. For $j\in[6]$, we let $\aaa_j=\aaa_0+\sum_{1\leq i\leq j}\uu_i\in A$. Since the distance between two red points in different staircases is strictly larger than $1$, the points $\aaa_0,\dots,\aaa_6$ must lie in the same red staircase.

Observe that for any $j\in[6]$, the $x$-coordinate of $\aaa_j$ exceeds that of $\aaa_{j-1}$ by at least $1/\sqrt{2}$. Hence, the $x$-coordinate of $\aaa_6$ exceeds that of $\aaa_0$ by at least $6/\sqrt{2}>\lambda+1$. Consequently, there exist points $\aaa_i$ and $\aaa_j$ in the horizontal stairs, $i<j$, such that $\aaa_j$ lies one level below $\aaa_i$, see~\Cref{stair_coloring} for an illustration. In particular, the difference between the $y$-coordinates of $\aaa_i$ and $\aaa_j$ is strictly larger than $\lambda-1=3/\sqrt{2}$. We choose such $i$ and $j$ to minimize $j-i$. In particular, $\aaa_{i+1}$ must lie in the vertical stair with width $1 < 2 \cdot 1/\sqrt{2}$, and thus $j-(i+1)\le 2$, i.e., $j \le i+ 3$. Since each of $\uu_1, \dots, \uu_6$ makes an angle of at most $\pi/4$ with the vector $(1,0)$, every two consecutive points in $\aaa_0,\dots,\aaa_6$ differ in the $y$-coordinate by at most $1/\sqrt{2}$. It then follows that the $y$-coordinates of $\aaa_i$ and $\aaa_j$ differ by at most $3/\sqrt{2}$, a contradiction.
\end{proof}

\begin{proof}[Proof of \Cref{leftovers}~(1)]
We shall give two constructions of $4$-colorings of $\R^2$ with no monochromatic unit-copy of $C_4$, i.e., showing that  $\chigraph{\R^2}{C_4}\le 4$. For the first one, we take any $a$ satisfying $1/\sqrt{3}<a<1/\sqrt{2}$, e.g., $a=2/3$, and consider a tiling of the plane by translates of a regular hexagon of side length $a$. Color these hexagons in 4 colors as in~\Cref{4coloring}~(1), where any point on the boundary of two or three tiles gets a color of one of these tiles arbitrarily. For the second construction, we consider a tiling of the plane by translates of a square of side length $1$ including the left and bottom side of its boundary and excluding other points on the boundary, i.e., the tile is $[0,1)\times [0,1)$. Color these tiles in 4 colors as in~\Cref{4coloring}~(2).

\begin{figure}[h!]
\begin{subfigure}
\centering
\begin{circuitikz}[scale=5/12]

\definecolor{color1}{HTML}{B2B2FF}
\definecolor{color4}{HTML}{FFCCCC}
\definecolor{color3}{HTML}{D25353}
\definecolor{color2}{HTML}{4747AF}

\def\hexsize{1} 
\def\dx{3*\hexsize/2} 
\def\dy{\hexsize*sqrt(3)} 

\def\nx{8}
\def\ny{6}

\foreach \i in {0,...,\nx} {
    \foreach \j in {0,...,\ny} {
        \pgfmathsetmacro{\x}{\i*\dx}
        \pgfmathsetmacro{\y}{\j*\dy + mod(\i,2)*\dy/2}
        
        \pgfmathtruncatemacro{\icol}{mod(\i,4)}
        \pgfmathtruncatemacro{\jcol}{mod(\j,2)}
        \ifcase\icol
            \ifcase\jcol \def\hexcolor{color1}\else \def\hexcolor{color2}\fi
        
        \or
            \ifcase\jcol \def\hexcolor{color3}\else \def\hexcolor{color4}\fi
        
        \or
            \ifcase\jcol \def\hexcolor{color2}\else \def\hexcolor{color1}\fi
        
        \or
            \ifcase\jcol \def\hexcolor{color4}\else \def\hexcolor{color3}\fi
        \fi
        
        \fill[fill=\hexcolor, draw=none] 
            (\x+\hexsize, \y) 
            \foreach \angle in {60,120,...,360} { -- ++(\angle:\hexsize) } -- cycle;
    }
}
\node at (14.3,0.5) {(1)};
\end{circuitikz}
\end{subfigure}
\hfill
\begin{subfigure}
\centering
\begin{circuitikz}[scale=2/3]

\definecolor{color1}{HTML}{B2B2FF}
\definecolor{color4}{HTML}{FFCCCC}
\definecolor{color3}{HTML}{D25353}
\definecolor{color2}{HTML}{4747AF}

\def\squaresize{1}

\def\nx{8}
\def\ny{6}

\foreach \i in {0,...,\nx} {
    \foreach \j in {0,...,\ny} {

        \pgfmathtruncatemacro{\icol}{mod(\i,2)}
        \pgfmathtruncatemacro{\jcol}{mod(\j,2)}
        \ifcase\icol
            \ifcase\jcol \def\squarecolor{color1}\else \def\squarecolor{color2}\fi
        \or
            \ifcase\jcol \def\squarecolor{color3}\else \def\squarecolor{color4}\fi
        \fi

        \pgfmathsetmacro{\x}{\i*\squaresize}
        \pgfmathsetmacro{\y}{\j*\squaresize}

        \fill[fill=\squarecolor, draw=none] (\x,\y) rectangle ++(\squaresize,\squaresize);
    }
}
\node at (9.6,0.3) {(2)};
\end{circuitikz}
\end{subfigure}
\caption{\!(1) a $4$-coloring of $\R^2$ via hexagonal tiling, \!(2) a $4$-coloring of $\R^2$ via square tiling.}
\label{4coloring}
\end{figure}
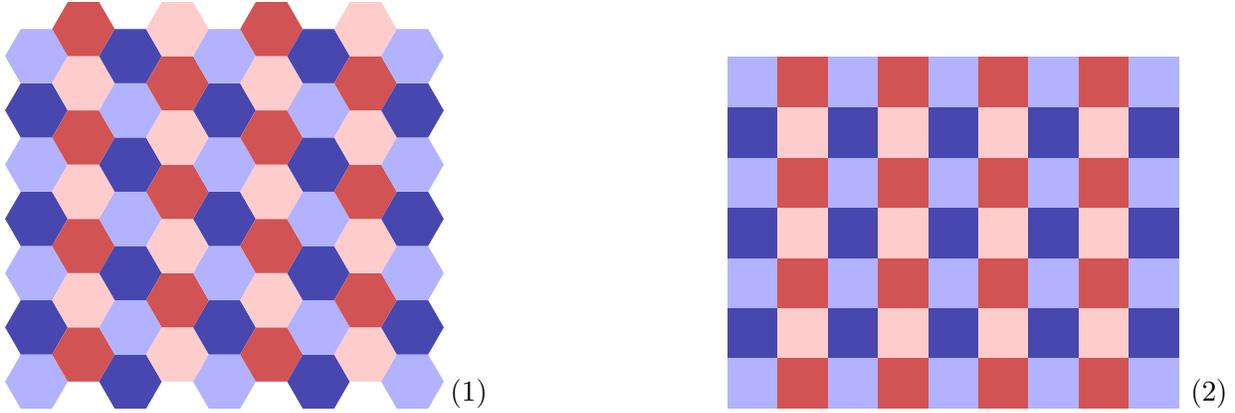

Note that a 4-point set $A$ is a unit-copy of $C_4$ if and only if $A$ is the vertex set of a rhombus of side length $1$. Now we show that no such $A$ is monochromatic under either of these colorings by considering the following two cases separately: (i) $A$ lies within one tile; (ii) $A$ is distributed between several tiles of the same color. 
On the one hand, observe that the distance between any two points in one tile of either coloring is strictly smaller than $\sqrt{2}$. However, the smallest diameter of a unit rhombus equals $\sqrt{2}$ and achieved only by squares. Hence, the case (i) is not possible. 
 On the other hand, observe that the distance between any two points of the same color from different tiles of either coloring is strictly greater than $1$. Hence, the case (ii) is not possible as well.
\end{proof}

\begin{proof}[Proof of \Cref{leftovers}~(2)]
Assume that $\chi(\R^2)=7$. Let $H$ be a finite unit-distance graph in the plane such that $\chi(H)=7$.  By~\Cref{product}, the Cartesian product  $H\square C_3$ is a unit-distance graph in the plane. Hence, to complete the proof, it is sufficient to show that any $2$-coloring of the vertices of $H\square C_3$ contains a monochromatic induced copy of $C_4$.
Suppose that the vertices of $H\square C_3$ are $2$-colored. For each $v \in V(H)$, consider the ``slice'' $\{v\}\square C_3$ of $ H\square C_3$. By the pigeonhole principle, this slice contains two vertices of the same color. Depending on which of the $3$ pairs of the vertices in $\{v\}\square C_3$ is monochromatic and which color this pair has, we define an auxiliary $6$-coloring of the vertices of $H$.
When all three vertices of the copy of $C_3$ are of the same color, we choose the pair of vertices defining the auxiliary coloring arbitrarily.
 Since $\chi(H)=7$, we conclude that some adjacent in $H$ vertices $v_1,v_2 $ are monochromatic under the auxiliary coloring, and thus the initial coloring of $H\square C_3$ contains a monochromatic induced copy of $C_4$, as desired.
\end{proof}

\begin{proof}[Proof of \Cref{leftovers}~(3)]
Now, we shall prove that $\chigraphind{\R^3}{C_4}>2$.
Note that this would be immediate using the argument from the proof of \Cref{leftovers}~(2) provided that $\chi(\R^3)\ge7$, but this is not known.
Suppose that the points of $\R^3$ are colored red and blue. One of the classical results in the field \cite[Theorem~8]{EGMRSS1} states that we can find a monochromatic copy of any triangle. Let $\x,\y,\z$ be a monochromatic, say red, triple of points such that $\|\x-\y\|=\|\y-\z\|=1, \|\x-\z\|=\sqrt{2}$. Consider the circle $C$  in the plane orthogonal to the plane containing the triangle $\x\y\z$, centered at the midpoint of the segment $\x\z$, and passing through $\y$. For every point $\y'$ on this circle except  for $\y$ and two points $\y_1,\y_2$ at distance $1$ from $\y$, the set $\{\x,\y,\z,\y'\}$ is an induced  unit-copy of $C_4$. If none of these copies is red, then all the points on $C$ except at most $3$ of them are blue. In this case, $C$ contains 4 blue vertices of a unit square, as desired.
\end{proof}

\section{Proofs of results on growing dimensions and of the canonical Ramsey-type result}
\label{sec:growing-canonical}

In this section we prove  \Cref{thm:main-asymptotic} and \Cref{thm:canon}.

For \Cref{thm:main-asymptotic}~(1), we need the following result on boxes.
A $d$-dimensional {\it box} is a Cartesian product of $d$ two-element sets, i.e., a product $\{0, b_1\} \times \cdots \times \{0, b_d\}$ for some positive real numbers $b_1, \ldots, b_d$ that we call \emph{side lengths} of the box.

\begin{lemma}[{Frankl--R\"odl~\cite[Lemma~3.1]{FR90}}] \label{almost_regular}
   Let $m \ge 11$. If $\varepsilon=\varepsilon(m)>0$ is sufficiently small, then for all choices of $b_{i,j}$, $1\le i < j \le m$, satisfying $|b_{i,j}-1|\le \varepsilon$, one can find a $\binom{m}{2}$-dimensional box $B$ and $m$ of its vertices $\x^{(i)} \in B$, $1\le i \le m$, such that $\|\x^{(i)}-\x^{(j)}\|^2=b_{i,j}$ for all $1\le i < j \le m$. If, in addition, all $b_{i,j}$ are rational numbers, then one can find such a box $B$ whose squared side lengths are also rational numbers.
\end{lemma}

\begin{proof}[Proof of \Cref{thm:main-asymptotic}~(1)]
By adding isolated vertices in $H$ if needed, we assume without loss of generality that $V(H)=[m]$ and $m \ge 11$. Let $\varepsilon=\varepsilon(m)$ be from the statement of~\Cref{almost_regular}, and for all $1\le i < j \le m$, let $b_{i,j}=1$ if $\{i,j\} \in E(H)$ and $b_{i,j}=1+\varepsilon$ otherwise. \Cref{almost_regular} implies that for $d=\binom{m}{2}$, there is a $d$-dimensional box $B \subseteq \R^d$ and $m$ of its vertices $\x^{(i)} \in B$, $1\le i \le m$, such that $\|\x^{(i)}-\x^{(j)}\|^2=b_{i,j}$ for all $1\le i < j \le m$.
Note that the set $A\coloneq \{\x^{(i)}: 1\le i \le m\} \subseteq B$ is an induced unit-copy of $H$ by construction, and thus $\chigraphind{\R^n}{H} \ge \chi(\R^n,B)$ for all $n$. Recall from \Cref{sec:background} that $\chi(\R^n,B)$ growth exponentially with $n$, which completes the proof.

We remark that as $n \to \infty$, Frankl and R\"odl~\cite{FR90} showed that a uniform lower bound $\chi(\R^n,B) \ge (b_d+o(1))^n$ holds for all $d$-dimensional boxes for some $b_d>1$. In fact, their methods yields that $b_d=1+(2+o(1))^{-d}$ suffices, see~\cite[Theorem~4]{AS18}. For the upper bounds, see \cite{Pros18}.
\end{proof}

\begin{proof}[Proof of \Cref{thm:main-asymptotic}~(2)]
For an $m$-partite $H$, let its vertex set be the union of pairwise disjoint independent sets $A_1, \ldots, A_m$. The idea of the proof is that we can find a unit-copy of $H$ in $\R^n$ such that $A_i$ lie on circles of radius $1/\sqrt{2}$ centered at the origin and contained in pairwise orthogonal planes.  
Thus, we can bound $\chigraph{\R^n}{H}$ from below with the generalized chromatic number of the orthogonality graph.

We say that a graph $G$ is an \emph{orthogonality graph in $\R^d$} if its vertices correspond to unit vectors in $\R^d$ such that adjacent vertices correspond to orthogonal vectors. We say that a graph $G$ is a \emph{graph of orthogonal planes in $\R^d$}, if its vertices correspond to planes in $\R^d$ containing the origin such that adjacent vertices correspond to orthogonal planes.
For each unit vector $\mathbf{v} \in \R^d$ consider  the plane $P(\mathbf{v})$ in $\R^{2d}$ spanned by two vectors $(\mathbf{v},\mathbf{0}), (\mathbf{0},\mathbf{v}) \in \R^{2d}$. 
Note that if $G$ is an orthogonality graph in $\R^d$, then $G$ is a graph of orthogonal planes in $\R^{2d}$.  Indeed, 
 two unit vectors  $\mathbf{v_1}$ and $\mathbf{v_2}$ in $\R^d$ are orthogonal if and only if $P(\mathbf{v_1})$ is orthogonal to $P(\mathbf{v_2})$. We remark that a similar idea was used in a different context in~\cite{BD23}.
 For a graph $G$, let the \emph{generalized chromatic number} $\chigraph{G}{K_m}$ be the minimum number of colors needed to color the vertices of $G$ such no copy of $K_m$ in $G$ is monochromatic. %

Let $G$ be a graph of orthogonal planes in $\R^n$ and $r= \chigraph{G}{K_m}-1$. We shall show that any $r$-coloring of $\R^n$ contains a monochromatic unit-copy of $H$.  
Define an auxiliary $r$-coloring of  $G$ as follows. Given a vertex $\mathbf{v}$ of $G$, let $C(\mathbf{v})$ be a circle of radius $1/\sqrt{2}$ in the plane $P(\mathbf{v})$ centered at the origin. Consider an arbitrary arc on $C(\mathbf{v})$ of diameter less than 1, say, of angular measure $\pi/4$.
This arc contains infinitely many points of the same color. If there are several such colors, we pick one of them arbitrarily. Color $\mathbf{v}$ in this color. By the choice of $r$, there are monochromatic, say red, vertices $\mathbf{v}_1, \ldots, \mathbf{v}_m$ forming a clique in $G$. Consider an injective embedding $f_i$ of the independent set $A_i$ into red points on the arc on $C(\mathbf{v}_i)$, for each $i=1, \ldots, m$. Since the diameter of the arc is less than 1, the distance between $f_i(v)$ and $f_i(u)$ is not 1 for any distinct $u, v\in A_i$. Since the distance between any point on $C(\mathbf{v}_i)$ and any point on $C(\mathbf{v}_j)$ is $1$,  for $1\leq i<j\leq m$, this embedding gives a red unit-copy of $H$. Therefore, $\chigraph{\R^n}{H}>r$, as desired. 
Moreover, if $H$ is a complete $m$-partite graph, then this embedding gives an induced unit-copy of $H$, and thus $\chigraphind{\R^n}{H} > r$.

To complete the proof, it remains only to show that some orthogonality graphs have generalized chromatic number exponentially large in terms of the dimension.
When $m=2$, we use a results by Raigorodskii~\cite{Rai99} claiming that there is a sequence $G_d$ of orthogonality graphs in $\R^d$ such that $\chi(G_d) \ge (2/\sqrt3+o(1))^d$ as $d \to \infty$.
For general $m$, we use a result by Frankl and R\"odl~\cite[Theorem~1.9]{FR87} implying that  for each $m \ge 2$, there exists $c'_m>1$ such that $\chigraph{J(4d,2d,d)}{K_m}>(c'_m+o(1))^{4d}$ as $d \to \infty$. 
Here, $J=J(n, k, t)$ is the \emph{generalized Johnson graph}  defined as a graph on the vertex set $\left\{\vv\in\{0,1\}^{n}:\, \|\vv\|^2=k\right\}$, where $\uu \vv$ is an edge if and only if $\uu$ and $\vv$ have precisely $t$ common ones.
One can see that $J(4d,2d,d)$ is an orthogonality graph in $\R^{4d}$ by replacing each 0 coordinate with~$-1$.
\end{proof}

We remark that the constants $c_m'$ produced by this method converge to $1$ very fast, namely, $c_m'=1+(2+o(1))^{-2^m}$ as $m \to \infty$, see~\cite[Lemma~3]{Sag18}.

\begin{proof}[Proof of \Cref{thm:main-asymptotic}~(3)]
For a graph $H$, define $\rho(H)$ as the infimum over all positive $\rho$ such that for a sufficiently large $n$, the $n$-dimensional ball $B_\rho^n \coloneq \{\x\in \R^n: \|\x\|\le \rho\}$ of radius $\rho$ contains a unit-copy of $H$ (not necessarily induced). Prosanov~\cite[Theorem~1]{Pros18} essentially proved that $\chigraph{\R^n}{H} \le (1+1/\rho(H) + o(1))^n$ for  any  graph $H$ as $n \to \infty$. Note that this implies the upper bound $\chigraph{\R^n}{H}\le (3+o(1))^n$ mentioned in \Cref{sec:basic-observations} via a straightforward inequality $\rho(H)\ge 1/2$, which is actually tight for all bipartite graphs.
We shall argue that $\rho(C_{2\ell+1}) = r_\ell \coloneq1/(2\cos(\frac{\pi}{4\ell+2}))$ for all $\ell \ge 1$.

\noindent
\textbf{Note added.} Mathieu Rundstrom pointed out to us that $2\rho^2(H) = 1-1/\theta(\bar H)$, where $\bar H$ is the complement of $H$ and $\theta(\cdot)$ is the Lov\'asz theta function, see \cite[p.~23]{Lov} or \cite[Theorem~3.1]{CST} for an explicit statement. Thus the value of $\rho(C_{2\ell+1})$ can be deduced from the existing literature. Nevertheless, we keep its short computation below for completeness. 

For the lower bound $\rho(C_{2\ell+1}) \ge r_\ell$, assume that for some $n$, there exist $2\ell+1$ points $\y_0, \y_1, \dots, \y_{2\ell}$ in $\R^n$ such that $\|\y_i\|<r_\ell<1/\sqrt2$ and $\|\y_{i}-\y_{i+1}\|=1$ for all $0\le i \le 2\ell$, where the indices are treated modulo $2\ell+1$. Observe that for all $0\le i \le 2\ell$, we have $\y_i\neq \mathbf{0}$ since otherwise $1=\|\y_{i}-\y_{i+1}\|=\|\y_{i+1}\|< r_\ell$, a contradiction. Given $0\le i \le 2\ell$, let $\theta_i=\angle\y_{i}\mathbf{0}\y_{i+1}$, where out of 2 possible values, we choose the one satisfying $0 \le \theta_i < \pi$.  The law of cosines implies that $1=\|\y_{i}-\y_{i+1}\|^2=\|\y_{i}\|^2+\|\y_{i+1}\|^2-2\|\y_{i}\|\|\y_{i+1}\|\cos\theta_i$. Therefore, $\cos\theta_i = -\frac{1-\|\y_{i}\|^2-\|\y_{i+1}\|^2}{2\|\y_{i}\|\|\y_{i+1}\|} < -\frac{1-2r_\ell^2}{2r_\ell^2} = -\cos(\frac{\pi}{2\ell+1}),$ and thus $\theta_i> \pi-\frac{\pi}{2\ell+1}$. Now the triangle inequality for angles  yields that $\angle\y_{i}\mathbf{0}\y_{i+2} \le \angle\y_{i}\mathbf{0}(-\y_{i+1})+\angle(-\y_{i+1})\mathbf{0}\y_{i+2}=(\pi -\theta_i)+(\pi-\theta_{i+1})< \frac{2\pi}{2\ell+1}$ for all $0\le i \le 2\ell-2$, and thus $\angle\y_{0}\mathbf{0}\y_{2\ell} \le \angle\y_{0}\mathbf{0}\y_{2}+\angle\y_{2}\mathbf{0}\y_{4}+\dots+\angle\y_{2\ell-2}\mathbf{0}\y_{2\ell} <\frac{2\ell\pi}{2\ell+1}$. On the other hand, $\angle\y_{0}\mathbf{0}\y_{2\ell} = \theta_{2\ell}>\pi-\frac{\pi}{2\ell+1}$ as shown above, a contradiction.

Even though the upper bound $\rho(C_{2\ell+1}) \le r_\ell$ is not needed for the proof of \Cref{thm:main-asymptotic}, we provide it here for completeness. Construct an explicit embedding of $C_{2\ell+1}$ onto a circle as a regular star polygon with the Schl\"afli symbol $\{2\ell+1/\ell\}$. In other words, let $\x_0,\x_1,\dots,\x_{2\ell}$ be the vertices of the regular $(2\ell+1)$-gon on the circle of radius $r_\ell$ centered at the origin $\mathbf{0} \in \R^2$. Connect $\x_i$ to $\x_{i+\ell}$ by a segment for all $0\le i \le 2\ell$, where the indices are treated modulo $2\ell+1$. We claim that all these segments are of unit length. Indeed, due to the symmetry of this construction, it is sufficient to check that $\|\x_0-\x_{\ell}\|= \|\x_0-\x_{\ell+1}\|=1$. Observe that $\alpha_\ell \coloneq \angle\x_\ell\mathbf{0}\x_{\ell+1}=\frac{2\pi}{2\ell+1}$, and thus $\|\x_\ell-\x_{\ell+1}\|=2r_\ell\sin(\frac{\alpha_\ell}{2}) = 2\sin(\frac{\pi}{4\ell+2})$. Besides, we have $\beta_\ell \coloneq\angle\x_\ell\x_0\x_{\ell+1}=\frac{\pi}{2\ell+1}$, and thus $\|\x_0-\x_{\ell}\|= \|\x_0-\x_{\ell+1}\| = \|\x_\ell-\x_{\ell+1}\|/2\sin({\frac{\beta_\ell}{2}})=1$, as desired.

To complete the proof, note that if $H$ is not bipartite, then $H$ contains an odd cycle $C_{2\ell+1}$ for some $\ell\ge 1$, and thus $\chigraph{\R^n}{H} \le \chigraph{\R^n}{C_{2\ell+1}}\le (3-\varepsilon_\ell+o(1))^n$, where $\varepsilon_\ell = 2-2\cos(\frac{\pi}{4\ell+2})>0$. We remark that $\varepsilon_\ell = (\frac{\pi}{4\ell})^2+O(\frac{1}{\ell^3})$ as $\ell \to \infty$.
\end{proof}

\begin{proof}[Proof of \Cref{thm:canon}]
For a graph $H$, we use \Cref{almost_regular} as in the proof of \Cref{thm:main-asymptotic}~(1) above but with a rational $\varepsilon$ to find a $d$-dimensional box $B$ with side length $b_1,\dots,b_d$ such that all $b_i^2$ are rational numbers and $B$ contains an induced unit-copy of $H$. Let $p_1,\dots ,p_d, q$ be positive integers such that $b_i^2 = p_i/q$ for each $i\in [d]$. Observe that for $p\coloneqq p_1+\dots+p_n$, the vertex set $Q\coloneq \{0, 1/\sqrt{q}\}^p$ of the $p$-dimensional hypercube of side length $1/\sqrt{q}$ contains a copy of $B$ and thus contains an induced unit-copy of $H$ as well. As mentioned in \Cref{sec:growing}, Geh\'er, Sagdeev, and T\'oth~\cite[Theorem~4]{GST} showed that there is a sufficiently large $n$ such that every coloring of $\R^n$ contains a copy $Q'$ of $Q$ that is either monochromatic or rainbow. Note that an induced unit-copy of $H$ inside $Q'$ is either monochromatic or rainbow as well.
\end{proof}

\section{Concluding remarks and open problems}\label{sec:conclusions}
\addtocontents{toc}{\protect\setcounter{tocdepth}{1}}

We initiated the study of a new Ramsey-type questions in Euclidean spaces. For a given graph $H$, we look for the largest integer $r$ such that in any $r$-coloring of $\R^n$ there is a monochromatic unit-copy of $H$. Many questions are open.

\subsection{Planar case} \Cref{leftovers}~(2) implies that if $\chi(\R^2)=7$, then $\R^2 \xrightarrow{2} C_4$. It would be very interesting to obtain this conclusion unconditionally.

\begin{question} \label{Q1}
Is it true that $\R^2 \xrightarrow{2} C_4$? 
\end{question}

\noindent
\textbf{Note added.}  Moore and Sagdeev~\cite{MS} answered \Cref{Q1} in the positive.

We know more about larger even cycles in the plane. \Cref{thm:main}~(2) implies that for $\ell=4$ and any integer $\ell \ge 6$, we have
$\R^2 \xrightarrow{r} C_{2\ell}$, where $ 5\leq r+1 =\chi(\R^2)\leq 7$. 
In particular, $\R^2 \xrightarrow{4} C_{2\ell}$. For $\ell = 3,5$, \Cref{thm:main-layer}~(1) implies that $\R^2 \xrightarrow{2} C_{2\ell}$, but the number of colors here is likely not the best possible. In particular, we cannot even show that $\R^2 \not\xrightarrow{6} C_6$ or $\R^2 \not\xrightarrow{6} C_{10}$.

\begin{question} \label{Q_C6C10}
Is it true that $\R^2 \xrightarrow{3} C_6$ and $\R^2 \xrightarrow{3} C_{10}$? 
\end{question}

\noindent
\textbf{Note added.} Alm\'asi, Axenovich, and Sagdeev~\cite[Corollary~1]{AlAxSag} answered \Cref{Q_C6C10} in a strong form by showing that $\chigraph{\R^2}{C_{6}} = \chigraph{\R^2}{C_{10}} = \chi(\R^2)$.

As for the odd cycles, \Cref{thm:main}~(3) shows that for all sufficiently large $\ell$, we have  $\R^2 \xrightarrow{r} C_{2\ell+1}$, where 
$r=\lceil \chi(\R^n)/ 2 \rceil-1$. 
Since $ 5\leq \chi(\R^2)\leq 7$, in particular, there exists $\ell_0$ such that $\R^2 \xrightarrow{2} C_{2\ell+1}$ for any integer $\ell\ge \ell_0$. To find some specific value for $\ell_0$,  we consider two special families of 5-chromatic unit-distance graphs in the plane, namely Heule graphs and Parts graphs, see~\cite{Heule, Parts}.
A computer verification shows that in every $2$-coloring of the vertices of \texttt{HeuleGraph610}, there is a monochromatic odd cycle of length at most 17.
Note that smaller Heule graphs and all Parts graphs do not have this property.
In particular this shows that 
$\R^2 \xrightarrow{2} \{C_3, C_5, \ldots, C_{17}\}$. Following the proof of \Cref{thm:main}~(3), it implies that $\R^2 \xrightarrow{2} C_{\ell}$ for any odd $\ell\geq 27$. On the other hand, recall from \Cref{sec:1.1} that $\R^2 \not\xrightarrow{2} C_3$. We do not know whether $\R^2 \not\xrightarrow{2} C_5$. 

\begin{question} \label{Q_odd}
What is the smallest  integer $\ell$,  $2\le \ell\le 13$ such that $\R^2 \xrightarrow{2} C_{2\ell+1}$? Is $\ell=2$?
\end{question}

Going back to bipartite graphs, \Cref{thm:main-layer}~(2) states that $\chigraph{\R^2}{Q_{d}}=2$ for $d=11$. We believe that a more careful analysis of staircase-like two-colorings can reduce the dimension perhaps to $d$ equal to $5$ or $6$, but probably not to $d=3$. 

\begin{question} \label{Q_bipartite}
What is the smallest order of a bipartite graph $H$ such that $\chigraph{\R^2}{H}=2$?
\end{question}

Some of the two-colorings of the plane allow to avoid monochromatic graphs whose every unit-copy is relatively rigid, i.e., whose unit-copies in the plane are almost isometric, or have relatively rigid subgraphs. For example, $K_3$ is rigid in this sense and the coloring by strips employs this property to avoid the monochromatic unit-copy of $K_3$. We observed by applying the pigeonhole principle to the edge-defining directions in a large hypercube, that it contains a sufficiently large subgraph whose vertices are close to points uniformly spaced on a straight line. We call it a ``pencil-like'' subgraph. That in turn allowed us to use a staircase-like two-coloring of the plane avoiding a monochromatic unit copy of a large hypercube.

\begin{question}
    Which triangle-free graphs $H$ are ``rigid'' as unit-distance graphs in the plane? Is it true that  $\R^2 \not\xrightarrow{2} H$ for such graphs?
    Which bipartite graphs, other than hypercubes, contain a large pencil-like subgraph  in every unit-copy in the plane? 
\end{question}

Note the grid $P_n\square P_n$ does not satisfy this condition: it can be embedded into three arbitrarily small discs centered at a copy of $T\coloneq\{0,1,2\} \subseteq \R$. Hence, the result of Currier, Moore, and Yip~\cite{CMY24} that $\R^2 \xrightarrow{2} T$ motivates the following strengthening of \Cref{Q1}.

\begin{question} \label{Q_Pn}
Is it true that $\R^2 \xrightarrow{2} P_n\square P_n$ for every $n \in \N$?
\end{question}

\subsection{Growing dimensions}
We provided a number of results for growing dimensions of the considered Euclidean space. These results have quite a different flavor and use different techniques compared to the fixed dimension ones.

Recall that $\chigraphind{\R^n}{H} \le \chigraph{\R^n}{H}$ for every graph $H$ and every $n\in \N$. In some degenerate cases, this inequality is not tight. Let $W_6$ be a wheel graph on 7 vertices, namely a graph obtained from $C_6$ by adding a new  vertex of degree 6. Let $H$ be a subgraph of $W_6$ obtained by removing one edge. There are two different options up to isomorphism, pick one arbitrarily. Then one can see that each unit-copy of $H$ in the plane is the vertex set of a regular unit hexagon plus its center. In particular, there are no induced unit-copies of $H$ in the plane and thus $\chigraphind{\R^2}{H}=1$, while $\chigraph{\R^2}{H}=2$. However, this effect seems to disappear as we increase the dimension.

\begin{question} \label{Q_ind}
    Is it true that for each graph $H$, we have $\chigraphind{\R^n}{H} = \chigraph{\R^n}{H}$ provided that $n$ is sufficiently large in terms of $H$?
\end{question}

\Cref{thm:main-asymptotic}~(2) provides a universal lower bound on $\chigraph{\R^n}{H}$ that holds for all $m$-partite graphs $H$. It is natural to ask for the smallest value of $\chigraph{\R^n}{H}$ among all $m$-partite graphs $H$. We suspect that this minimum may be equal to $\chigraph{\R^n}{K_m}$ provided that $n$ is sufficiently large. In particular, this would imply that $\chigraph{\R^n}{C_4} = \chi(\R^n)$ and $\smash{\chigraph{\R^n}{K_4^-}}= \chigraph{\R^n}{K_3}$ for all sufficiently large $n$, where $K_4^-$ is the graph obtained from $K_4$ by removing one edge. It seems that our ``Cartesian'' approach does not suffice to get these equalities.

\begin{question} \label{Q_partite}
    Is it true that $\chigraph{\R^n}{H}\ge \chigraph{\R^n}{K_m}$ for each $m \ge 2$ and each $m$-partite graph $H$ provided that $n$ is sufficiently large in terms of $H$? In particular, is it true that $\chigraph{\R^n}{H}= \chi(\R^n)$ for each bipartite graph $H$ provided that $n$ is sufficiently large in terms of $H$?
\end{question}

By considering a proper coloring of $\R^n$ with $\chi(\R^n)$ colors and then grouping these colors into $\lceil \chi(\R^n)/m \rceil$ classes, each containing at most $m$ colors, we see that $\chigraph{\R^n}{H}\le \lceil \chi(\R^n)/m \rceil$ for every $(m+1)$-chromatic graph $H$. \Cref{thm:main}~(1) and~(2) demonstrate that for some $2$-chromatic graphs, this upper bound is tight for all $n \ge 2$.~\Cref{thm:main}~(3) shows that for every $n\geq2$ there are non-bipartite graphs attaining the upper bound. However, it is not clear whether this can be tight for a fixed non-bipartite graph as $n$ grows, and \Cref{thm:main-asymptotic}~(3) suggests that $\chigraph{\R^n}{H}$ could be exponentially smaller than $\chi(\R^n)$ in this case.

\begin{question}
    Is there $m \ge 2$ and an $(m+1)$-chromatic graph $H$ such that $\chigraph{\R^n}{H}= \lceil \chi(\R^n)/m \rceil$ for all $n\ge 2$?
\end{question}

Observe that for all graphs $H_1,H_2$, both $H_1$ and $H_2$ are subgraphs of their disjoint union $H_1 \sqcup H_2$ and thus $\chigraph{\R^n}{H_1 \sqcup H_2} \le \min \{\chigraph{\R^n}{H_1}, \chigraph{\R^n}{H_2}\}$ for all $n$. In some special cases, e.g., if $H_1=H_2$, a simple application of the pigeonhole implies that this inequality is tight, but we do not know if this is always the case.

\begin{question}
    Is it true that $\chigraph{\R^n}{H_1 \sqcup H_2} = \min \{\chigraph{\R^n}{H_1}, \chigraph{\R^n}{H_2}\}$ at least for large enough $n$?
\end{question}

For a graph $H$, let $d(H)$ be the minimum $d$ such that there is a $d$-dimensional box that contains an induced unit-copy of $H$. This parameter controls the lower bound on $\chigraphind{\R^n}{H}$ that our proof of \Cref{thm:main-asymptotic}~(1) leads to. \Cref{almost_regular} implies that $d(H) \le \binom{m}{2}$ for each graph $H$ on $m\ge 11$ vertices. We can use an embedding constructed for a different problem in~\cite{FPRR17} to obtain another upper bound $d(H)\le m(\overline{e}+1)$ for all graphs on $m$ vertices with $\binom{m}{2}-\overline{e}$ edges, which is stronger than the former bound for ``almost complete'' graphs. None of the upper bounds are tight in general, so it would be interesting to improve upon them and to study the relation between $d(H)$ and other standard graph parameters.

Recall from \Cref{sec:growing-canonical} that for a graph $H$, we define $\rho(H)$ as the infimum over all positive $\rho$ such that for a sufficiently large $n$, the $n$-dimensional ball of radius $\rho$ contains a unit-copy of $H$ (not necessarily induced). This parameter controls the upper bound $\chigraph{\R^n}{H}\le (1+1/\rho(H)+o(1))^n$ our proof of \Cref{thm:main-asymptotic}~(3) leads to. It is not hard to see that $\rho(H)=1/2$ for all bipartite graphs and that $\rho(K_m)=\sqrt{(m-1)/(2m)}$ for cliques. Hence, for each graph on $m$ vertices, we have $\rho(H)\le \rho(K_m)<1/\sqrt{2}$, and thus this method cannot produce an upper bound on $\chigraph{\R^n}{H}$ stronger than $(1+\sqrt{2}+o(1))^n$. We have shown that $\rho(C_{2\ell+1}) = 1/(2\cos(\frac{\pi}{4\ell+2}))$ for odd cycles. It would be interesting to compute this parameter for other natural classes of graphs or at least provide an algorithm for that.

\begin{question}
    Is $\rho(H)$ computable as a function of $H$?
\end{question}

\noindent
\textbf{Note added.} Mathieu Rundstrom pointed out to us that $2\rho^2(H) = 1-1/\theta(\bar H)$, where $\bar H$ is the complement of $H$ and $\theta(\cdot)$ is the Lov\'asz theta function, see \cite[p.~23]{Lov} or \cite[Theorem~3.1]{CST} for an explicit statement. In particular, the value of $\rho(H)$ can be computed in polynomial time.

\subsection{Infinite graphs}

All results we stated earlier concern only finite graphs. Nevertheless, it is reasonable to study the values of $\chigraph{\R^n}{H}$ and $\chigraphind{\R^n}{H}$ for infinite graphs as well. We claim that our proof of \Cref{thm:main-asymptotic}~(2) works in this case verbatim and yields the following.

\begin{theorem} \label{thm:inf_asymp}
    For each $m>2$, there exists $c_m>1$ such that $\chigraph{\R^n}{H} \ge (c_m+o(1))^n$ for every $m$-partite graph $H$ whose vertex set has cardinality at most that of the continuum. Moreover, for $m=2$, one can take $c_2= (4/3)^{1/4} \sim 1.074\dots$
\end{theorem}

Observe that \Cref{thm:inf_asymp} is in a sharp contrast with a classical result that only finite sets can be Ramsey. In fact, it follows from \cite[Theorem~19]{EGMRSS2} that $\chi(\R^n,A)=2$ for every infinite $A \subseteq \R^n$.

We also note that if a graph $H$ does not satisfy the conditions of \Cref{thm:inf_asymp}, namely if it has an infinite chromatic number or the cardinality of its vertex set is greater than that of the continuum, then it is not difficult to show that $\R^n$ does not contain a unit-copy of $H$ for any $n \in \N$, and thus $\chigraph{\R^n}{H}=1$ holds trivially.

It would be interesting to find an analogue of \Cref{thm:inf_asymp} for the induced variant. Note that our proof of \Cref{thm:main-asymptotic}~(1) does not work for infinite graphs.

\begin{question} \label{Q_inf_induced}
    What is the necessary and sufficient condition for an infinite graph $H$ to have an induced unit-copy in $\R^n$ for some $n$? Does $\chigraphind{\R^n}{H}$ grow with $n$ for such graphs?
\end{question}

Our only result on $\chigraphind{\R^n}{H}$ for infinite graphs is the following.

\begin{proposition}\label{inf_tree}
    Let $K_{1, \aleph_0}$ be a star with countably infinite number of leaves. Then $\R^2 \ind{2} K_{1, \aleph_0}$.
\end{proposition}

\begin{proof}[Proof of~\Cref{inf_tree}]
Assume the contrary, namely that there exists a red-blue coloring of the plane such that no induced copy of $K_{1, \aleph_0}$ is monochromatic. In particular, this implies that any unit circle contains only finitely many points of the same color as the center of the circle. 
Without loss of generality, suppose that the origin $\mathbf{0}$ is colored red. For a point $\z \in \R^2$, let $C(\z)$ be the unit circle centered at $\z$. Let $\x_0 \in C(\mathbf{0})$ be blue and $\x_1, \x_2,\dots$ be an infinite sequence of blue points on $C(\mathbf{0})\sm C(\x_0)$. Let $Y = \{\y \in C(\x_0) : \|\y\| < 1\}$. Observe that for each $i\ge 1$, there are only finitely many blue points on $C(\x_i)$, and each of them lies on at most 2 circles $C(\y), \y \in Y$. Hence, the set $Y'$ of points $\y \in Y$ such that $C(\y)$ intersects some $C(\x_i)$, $i \ge 1$, by a blue point is countable.
Recall that $Y$ contains only finitely many blue points. In particular, there exists a red point $\y_0 \in Y\sm Y'$. Observe that for all $i \ge 1$, we have $\|\y_0-\x_i\|<2$, and thus $C(\y_0)$ intersects $C(\x_i)$ at two points, both of which are red by construction. It remains to note that each point on $C(\y_0)$ belongs to at most two different circles $C(\x_i)$, $i \ge 1$, because $\mathbf{0} \notin C(\y_0)$. Therefore, the circle $C(\y_0)$ contains infinitely many red points, while its center $\y_0$ is also red, a contradiction.
\end{proof}

We suspect that a much stronger form of \Cref{inf_tree} might hold.

\begin{question}
    Is it true that $\chigraphind{\R^n}{F}\!=\!\chigraph{\R^n}{F}\!=\!\chi(\R^n)$ for any countable forest $F$ and any $n \ge 2$?
\end{question}

\subsection{Other norms} Analogues of the Euclidean Ramsey theory questions for other normed spaces have been studied as well, see for example~\cite{davies2025spacetime, FGST, FKS, Geh23}. We remark that the methods from \Cref{sec:growing-canonical} can be directly applied for other norms. However, the applications of our ``Cartesian'' approach are limited, since even for other $\ell_p$-norms, a Cartesian power of a unit-distance graph in $\R^n$ is not necessarily a unit-distance graph in the same dimension.

\subsection{Cartesian powers} Informally speaking, \Cref{thm:Turan} implies that if $G$ is a graph, $r$ is a ``small'' positive integer, and $H$ is a graph with some Tur\'an-type density equal to zero, then $G^{\square N} \xrightarrow{r} H$ for a sufficiently large integer $N$. However, \Cref{lem:layered} shows that even in the case $G=K_3$ and $r=2$, the density condition can be substantially relaxed to $H$ being layered. We wonder if this condition can be relaxed even further.

\begin{question} \label{Q14}
    For which graphs $H$ there exists a sufficiently large integer $N$ such that $K_3^{\square N} \xrightarrow{2}H$? 
\end{question}

We remark that $K_3^{\square N} \not\xrightarrow{2}H$ when $H=C_4$.
Indeed, for any $N\in \N$, we identify the vertices of $K_3^{\square N}$ with $[3]^n$. Then we color each $(x_1,\dots,x_N)\in [3]^N$ red if $x_1+\dots+x_N = 0 \pmod{3}$ and blue otherwise. One can check that this coloring contains no monochromatic copies of $C_4$. 

\noindent
\textbf{Note added.} Alm\'asi, Axenovich, and Sagdeev~\cite[Proposition~1]{AlAxSag} answered \Cref{Q14} by proving that the condition in \Cref{lem:layered} is best possible. Namely, they showed that if $K_3^{\square N} \xrightarrow{2}H$ for some $N$, then $H$ is layered.
\\

\noindent
{\bf Acknowledgements.} We thank Imre Leader for discussions on the Block Sets Conjecture, see~\cite{ILW}, and for spotting a mistake in one of our earlier arguments. We are grateful to J\'anos Pach for conversations on known Ramsey-type results for Cartesian powers of graphs. We also thank Kenneth Moore for comments on his online catalog of dense unit-distance graphs in the plane~\cite{Moore}. We thank Mathieu Rundstrom for bringing the references~\cite{CST, Lov} to our attention. We thank an anonymous reviewer for bringing the reference~\cite{HY} to our attention.

\noindent
\textbf{AI disclosure.} No AI assistance was used in obtaining the results of this paper or in its writing.

\addtocontents{toc}{\protect\setcounter{tocdepth}{2}}

\appendix

\section{Induced zero Tur\'an density } \label{app:hypercube-Turan}

In this section, we identify the vertices of a hypercube $Q_n$ with the subsets of the ground set $[n]$, where the two vertices are adjacent if and only if the respective sets have symmetric difference of size 1.
Let $L_k=L_k(Q)$ be the $k$th edge layer of a hypercube $Q=Q_n$ induced by  $\binom{[n]}k\cup \binom{[n]}{k-1}$. Recall that a graph is {\it layered} if it is a subgraph of an edge layer of a hypercube. We visualize the vertex layers of the hypercube arranged on different  horizontal lines such that the $i$th vertex layer line is above the $j$th vertex layer line for $i>j$.  
We call a star with $k$ edges and center in a higher vertex layer, than the vertex layer containing the leaves a {\it down $k$-star}.
We see that each vertex in the top vertex layer of $L_k$ is incident to exactly $k$ edges of $Q_n$ and corresponds to a down $k$-star.

For a set $A$ of vertices from the $k$th vertex layer of a hypercube, we associate vertices $A$ with edges of a  $k$-uniform hypergraph $\cH_A$.
We say that a graph $G$ has a {\it partite representation} if there are integers $k$ and $n$ and a set $A$ of vertices from the $k$th vertex layer of $Q_n$ such that 
$G$ is isomorphic to a subgraph $G'$ of a union of down $k$-stars centered at vertices from $A$ and such that the respective hypergraph $\cH_A$ is $k$-partite. If $G'$  is induced in $Q_n$, we say that $G$ has an  {\it induced partite representation}.

Conlon \cite{C} proved that if a graph $G$ has partite representation, it has zero hypercube Tur\'an density, that is $\ex(Q_n, G)=o(|E(Q_n)|)$.
We note here that the same statement holds in an induced setting and follows from the proof by Conlon.  We state and outline its proof here for completeness.

\begin{theorem}\label{thm:a1}
Let $G$ have an induced $k$-partite representation in the hypercube. 
Then a largest subgraph of $Q_n$ with no induced copy of $G$ has $o(|E(Q_n)|$ edges.
\end{theorem}

\begin{proof}[Proof of~\Cref{thm:a1}]
Let $\gamma>0$ and $n\in\mathbb{N}$ be sufficiently large.
Let $F$ be a subgraph of $Q_n$ on $\gamma |E(Q_n)|$ edges. 
By a standard argument and pigeonhole principle, see for example Conlon \cite{C}, we can assume that there is an edge layer $L_m$ of $Q_n$, $n/4\leq m \leq 3n/4$, such that $F$ contains at least $\gamma' |E(L_m)|$ edges of $L_m$, for $\gamma' > \gamma/4$.
Let $X= \binom{[n]}{m}$, the set of vertices of the top vertex layer of $L_m$.

Consider all sub-hypercubes of $Q_n$ with the ``lowest" vertex in $\binom{[n]}{m-k}$ and the ``highest" vertex $[n]$, that is, all $n'$-dimensional sub-hypercubes of $Q_n$ containing $[n]$ for $n'=n-m+k$. Denote the set of these sub-hypercubes $\cQ$. 
We see that the edge layer $L_m$ of $Q_n$ is the union of the  $k$th edge layers of the hypercubes from $\cQ$.
We shall argue that for a constant $\mu= (\gamma'/k)^k$ and  some $Q\in \cQ$, 
there is at least a $\mu$-proportion of down $k$-stars with edges from $F$ in the $k$th edge layer of $Q$.
Note that each down $k$-star from $L_m$  is a down $k$-star in exactly one $Q\in \cQ$.

Counting down $k$-stars in $L_m\cap F$ and applying the pigeonhole principle, we see that there is a $Q\in \cQ$ with at least $x$ down $k$-stars in its $k$th edge layer,  where
$$x= \binom{\gamma'm}{k} \binom{n}{m} / \binom{n}{m-k}\geq (\gamma'/k)^k \binom{n-m+k}{k} = \mu \binom{n'}{k}.$$
Here, we used the fact that the average degree of vertices from $\binom{[n]}{m}$ in $L_m$ is at least $\gamma' m$.  Let $Q'$ be a copy of $Q$ on the ground set $[n']$ and let 
 $A$ be the set of centers of down $k$-stars  in the $k$th edge layer of $Q'$.
Let $\cH=\cH_A$ be the $k$-uniform hypergraph on vertex set $[n']$ associated with $A$. Since $n'$ is sufficienlty large,  by a result of Erd\H{o}s, $\cH$ contains a copy of the  $k$-partite hypergraph $\cH'$ giving induced representation of $G$.
Translating back to the vertices from $\binom{[n']}{k}$, let $A'$ be the set of vertices corresponding to the edges of $\cH'$.
Then the union $F'$ of down $k$-stars centered at $A'$ is an induced subgraph of $Q'$ that corresponds to an induced subgraph of $Q$ and thus of $Q_n$. Since $F'$ contains an induced copy of $G$, so does $F$.
 \end{proof}

The following three partite representation is given in \cite{A} without pointing out that it is in fact an induced representation. We include it here for completeness as well.

\begin{lemma}\label{lem:a2}
Every even cycle of length $8$ or at least $12$ has an induced  three-partite or two-partite representation. 
\end{lemma}

\begin{proof}[Proof of~\Cref{lem:a2}]
Let $k\in \mathbb N$, $k\geq 3$. We shall construct a cycle $C=C_{2\ell}$ for odd $\ell=2k+1$ with vertices in the third and second vertex layer of $Q_n$. We shall define $A_\ell$  to be the ordered set of vertices of $C$ in the third layer, $B_\ell$ to  be a respective set of vertices in the second layer, and let  $C$ be  built by taking vertices from $A_\ell$ and $B_\ell$ alternatingly in order.
Let 
\begin{eqnarray*}
A_{2k+1} & = & (ax_1y_1,    ax_2y_1,  ax_2y_2,  \ldots, ax_{k-1}y_{k-2},  ax_{k-1}y_{k-1}, ~~ bx_{k-1}y_{k-1},   bx_{k-1}y_{0},  bx_1y_0,  bx_1y_1),\\
B_{2k+1} & = & (ay_1,    ax_2,  ay_2,  \ldots, ax_{k-1}, ~~  x_{k-1}y_{k-1},   bx_{k-1},  by_0,   bx_1,  x_1y_1).
\end{eqnarray*}
We see that $C$ is induced by $A_\ell\cup B_\ell$. 
Moreover, the elements of $A_\ell$ correspond to the hyperedges of a $3$-partite hypergraph with parts $\{a,b\}$, $\{x_1, x_2, \ldots, x_{k-1}\}$, and $\{y_0, y_1, y_2, \ldots, y_{k-1}\}$.

For $C=C_{2\ell}$ with $\ell =2k$, $k\geq 2$, $k\in \mathbb{N}$,
Consider ordered sets $A_{2k}$ and $B_{2k}$  of sizes $2k$ each, in the second and first vertex layers of $Q_n$, respectively, and let vertices of $C$ alternate between $A_{2k}$ and $B_{2k}$. 
Let 
\begin{eqnarray*}
A_{2k} & = & (x_0x_1,    x_1x_2,  x_2x_3,  \ldots, x_{k-1}x_{0}), \\
B_{2k} & = & (x_0, x_1,  x_2,  \ldots, x_{k-1}).
\end{eqnarray*}
We see that $C$ is induced by the vertices from $A_\ell$ in the second vertex layer and all their neighbors in the first layer. Moreover, the elements of $A_\ell$ correspond to the edges of a bipartite graph.
\end{proof}

However, one can not extend the above lemma to all layered subgraphs of a hypercube.

\begin{lemma}\label{lem:a3}
There are layered graphs that are not induced subgraphs of any edge layer of any hypercube.
\end{lemma}

\begin{proof}[Proof of~\Cref{lem:a3}]
Let $t\geq 4$ and  $H(t)$ be a graph that is a union of $t-1$ paths of length $5$ that pairwise share only endpoints, called {\it poles}.
We claim that $H(t)$ is a layered graph and in any its embedding in an edge layer of a hypercube the poles are adjacent.
To see that $H(t)$ is layered, embed the poles in  a vertex $[t]$ and $[t+1]$ of $Q_n$ on the ground set $[2t+1]$, respectively. 
Let the $(j-1)$st path between the poles have internal vertices $[t]\cup \{t+j\},  [t]\sm\{j\}\cup \{t+j\}, [t]\sm\{j\}\cup\{t+j, t+1\}, [t+1]\sm\{j\}$, for $j=2, \ldots, t$.

Assume now there is an embedding of $H(t)$ into an edge layer of $Q_n$ with poles $X$ and $Y$, $|X|=|Y|-1$,  such that the poles are not adjacent, i.e., such that $X$ is not a subset of $Y$. 
Since $X$ and $Y$ are connected by a path of length $5$, the symmetric difference $X\Delta Y$ is of odd size at most $5$.

If $|X\Delta Y|=5$, then $|Y\sm X|= 3$ and each vertex adjacent to $Y$ on an $X,Y$-path is equal to $Y-\{y\}$, for some $y\in Y\sm X$. Thus there are at most three internally vertex-disjoint $X,Y$-paths of length $5$, a contradiction. 

If $|X\Delta Y|=3$, then let $Y\sm X= \{y_1, y_2\}$ and $X\sm Y=\{x\}$. Then any $X,Y$-path in $Q_n$ has edges with directions $y_1, y_2, x, j, j$, for $j\in X\cap Y$. Here, a direction of an edge is the element in the difference of its endpoints. The only way to have such directions in order on $X,Y$-path is  
$y_2, j, y_1, x, j$ or $y_1, j, y_2, x, j$. In particular, the first internal vertex on all these paths is either $X\cup \{y_1\}$ or $X\cup \{y_2\}$. Thus the number of such internally vertex-disjoint paths is at most two, a contradiction. 
\end{proof}

\vspace{-1em}\section{Explicit exponential lower bounds on $\chigraph{\R^n}{H}$ and $\chigraphind{\R^n}{H}$}\label{app:exponential-lb}
\addtocontents{toc}{\protect\setcounter{tocdepth}{1}}

In this section, we obtain lower bounds on $\chigraph{\R^n}{H}$ and $\chigraphind{\R^n}{H}$ that are asymptotically better than those presented in \Cref{sec:growing} for some special classes of graphs.
These classes include graphs with slightly superlinear extremal numbers, graphs with orthogonal tree representations, and $C_4$.

\subsection{Extremal numbers}  \label{sec:b1}

Raigorodskii~\cite{Rai00} constructed unit-distance graphs with small independence number. 
As a consequence, he obtained the best known lower bound~\eqref{chiRn} on $\chi(\R^n)$.
We use these graphs to obtain essentially the same lower bound on $\chigraph{\R^n}{H}$ for sparse graphs $H$ with ``almost linear'' extremal number. For a graph $H$ and a positive integer $m$, the \emph{extremal number} $\ex(m,H)$ is defined as the maximum number of edges in a graph on $m$ vertices containing no subgraph isomorphic to $H$. It is well-known that $\ex(m,H) = \Theta(m^2)$ for each non-bipartite graph $H$ and $\ex(m,H) = o(m^2)$ for each bipartite graph $H$ as $m$ grows, see the survey by F\"uredi and Simonovits~\cite{FS}.

\begin{theorem}\label{linear_turan}
    Let $\varepsilon$ be any number satisfying $0< \varepsilon < 0.2$ and $\psi \coloneq 1.239\ldots$ be from~\cite{Rai00}.  Then there is $\delta>0$ such that the following holds. Given a graph $H$, if there is $m_0 \in \N$ such that $\ex(m,H) < m^{1+\delta}$ for all $m \geq m_0$, then there is $n_0 \in \N$ such that $\chigraph{\R^n}{H} \ge (\psi-\varepsilon)^n$ for all $n \geq n_0$.
\end{theorem}
\begin{proof}[Proof of~\Cref{linear_turan}]
    Let $0<\varepsilon< 0.2$, fix a sufficiently small $\delta>0$ such that $\frac{\psi-\frac{\varepsilon}{2} }{\psi-\varepsilon} > (\frac{2.277}{\psi-\varepsilon})^\delta$. 
    Let $H$ be a graph such that $\ex(m,H) < m^{1+\delta}$ for all sufficiently large $m$. Let $n$ be a sufficiently large integer.
    We omit floors and ceiling when clear from context. Assume for the sake of contradiction that $\chigraph{\R^n}{H} < (\psi-\varepsilon)^n$. 
       
    Let $G_n$ be the unit-distance graph in $\R^n$  with 
    $2.275^n\leq |V(G_n)| \leq 2.276^n$ and $\alpha(G) < \frac{|V(G)|}{(\psi- \frac{\varepsilon}{2} )^n}$, 
    constructed by Raigorodskii \cite{Rai00}. Consider  a coloring of $V(G_n)$ in less than $(\psi-\varepsilon)^n$ colors such that no copy of $H$ in $G_n$ is monochromatic. By the pigeonhole principle, this coloring contains a monochromatic subset $W$ such that $|W|= \frac{|V(G_n)|}{(\psi-\varepsilon)^n}$. We see that 
    $ 2\alpha(G) \leq |W| \leq  ( \frac{2.276}{\psi-\varepsilon})^n$.  
    By a Tur\'an-type argument applied to the subgraph $F$ of $G_n$ induced by $W$
    we see that $|E(F)|\geq \frac{|W|^2}{(4\alpha(G_n))}$. Observe that

    \begin{align*}
        \frac{|E(F)|}{|W|^{1+\delta}} &\geq \frac{|W|^2}{4\alpha(G_n)|W|^{1+\delta}}  =\frac{1}{4}\left(\frac{|W|}{|V(G_n)|}\right)\left(\frac{|V(G_n)|}{\alpha(G_n)}\right)\frac{1}{|W|^{\delta}} \\
        &\geq \frac{1}{4}\left(\frac{1}{\psi-\varepsilon}\right)^n 
    \left(\frac{\psi-\frac{\varepsilon}{2}}{1}\right)^n 
    \left(\frac{\psi-\varepsilon}{2.276}\right)^{\delta n}
    \geq \left( \left( \frac{\psi-\frac{\varepsilon}{2}}{\psi-\varepsilon}\right) 
     \left(\frac{\psi-\varepsilon}{2.277}\right)^{\delta} \right)^n > 1.
    \end{align*}
    In particular, $|E(F)|>|W|^{1+\delta}$ and thus $F$ contains a copy of $H$, a contradiction. 
\end{proof}

\subsection{Orthogonal trees}

Given a tree $T$ with $d$ edges, an \emph{orthogonal copy of $T$} is the vertex set $V(T')$ of a copy $T'$ of $T$ embedded into $\R^d$ such that all edges of $T'$ are pairwise orthogonal. Kupavskii, Zakharov, and the third author showed in~\cite{KSZ22} that for every tree $T$ with $d$ edges and every orthogonal copy $A$ of $T$, we have $\chi(\R^n, A) \ge (\psi^{1/d}+o(1))^n$ as $n \to \infty$, where $\psi=1.239\ldots$ is from \Cref{sec:b1}.
This immediately gives the following.

\begin{theorem}
Let a graph $H$ have an induced unit-copy in an orthogonal copy $A$ of a tree in $\R^d$. Then 
$\chigraphind{\R^n}{H} \ge \chi(\R^n, A) \ge (\psi^{1/d}+o(1))^n,$ as $n$ grows, where $\psi \coloneq 1.239\ldots$ is from~\cite{Rai00}.
\end{theorem}

Note that orthogonal copies of trees can contain induced unit-copies of ``dense'' graphs. 
For instance, the  leaf set of a star $K_{1,d}$ embedded into $\R^{d}$ such that all its edges are pairwise orthogonal and of length $1/\sqrt{2}$ forms a unit-copy of $K_d$, and thus $\chigraphind{\R^n}{K_d} \ge (\psi^{1/d}+o(1))^n$, which is asymptotically the best known lower bound on this chromatic number. As another example, let $T$ be a tree obtained from a star $K_{1,d}$ by connecting one of its leaves to two new vertices. It is not hard to find an orthogonal copy of $T$ in $\R^{d+2}$ containing an induced unit-copy of $K_{d+1}^-$, the graph obtained from $K_{d+1}$ by removing one edge. Therefore, we have $\chigraphind{\R^n}{K_{d+1}^-} \ge (\psi^{1/(d+2)}+o(1))^n$. Note that for large $d$, the last lower bound is only slightly weaker that the lower bound on  $\chigraphind{\R^n}{K_d}$ mentioned above. We suspect that even the equality $\chigraphind{\R^n}{K_{d+1}^-} = \chigraphind{\R^n}{K_{d}}$ may hold for all fixed $d$ and all sufficiently large $n$, see Questions~\ref{Q_ind} and~\ref{Q_partite}.

However, not all graphs have such orthogonal tree representation, and it would be interesting to describe those that have it. We verified that all graphs on at most four vertices have such a representation, and there are precisely three graphs on five vertices that do not: (1) the \emph{gem}, i.e., the graph obtained from $P_4$ by adding a new vertex of degree $4$; (2) the \emph{house}, i.e., the union of $C_4$ and $C_3$ sharing $2$ vertices and $1$ edge; (3) the graph obtained from $K_4$ by adding a new vertex of degree $2$.

\begin{question}
    What is the necessary and sufficient condition for a graph to admit an induced unit-copy in an orthogonal copy of some tree? 
\end{question}

\noindent
\textbf{Note added.} It can be shown that a graph admits an induced unit-copy in an orthogonal copy of some tree if and only if this graph is an \emph{exact $k$-leaf power} for some $k \in \N$, see~\cite{BR}. See also~\cite{HY} for a closely related yet deferent problem on realizing a given distance matrix by a tree.

\subsection{Frankl--R\"odl product argument}

Frankl and R\"odl used the product argument~\cite[Theorem~2.2]{FR90} to prove that if $A$ is a vertex set of a non-degenerate simplex, then $\chi(\R^n,A)$ grows exponentially with $n$. Applied to our graph setting, this argument leads to explicit exponential lower bounds on $\chigraphind{\R^n}{H}$ for graphs of the form $H=H_1\square H_2$. In particular, we obtain the following improvement upon the bound from \Cref{thm:main-asymptotic}~(2) for $C_4=K_2\square K_2$.

\begin{theorem}\label{thm:b3}
We have $\chigraphind{\R^n}{C_4}\geq\left(1.0792...+o(1)\right)^n$ as $n \to \infty$.
\end{theorem}

\begin{proof}[Proof of~\Cref{thm:b3}]
We repeat the argument of Frankl and R\"odl with minor modifications.  Let $G$ and $H$ be two unit-distance graphs in $\R^n$. Let $\alphagraphind{G\square H}{C_4}$ be the largest order of an induced subgraph of $G\square H$ which contains no induced copy of $C_4$. By the pigeonhole principle, one needs at least $|V(G\square H)|/\alphagraphind{G\square H}{C_4}$
colors to color $V(G\square H)$ such that no induced copy of $C_4$ is monochromatic. Since $G\square H$ is a unit-distance graph in $\R^n$ by~\Cref{product}, it follows that 
\begin{equation}
\label{chromatic}
\chigraphind{\R^n}{C_4}\geq\frac{|V(G\square H)|}{\alphagraphind{G\square H}{C_4}}.
\end{equation}
Our proof consists of two steps. First we obtain an upper bound for $\alphagraphind{G\square H}{C_4}$. Then we optimize the resulting lower bound on $\chigraphind{\R^n}{C_4}$ by choosing $G$ and $H$ carefully.

Let $W$ be an induced subgraph of $G\square H$ of order $\alphagraphind{G\square H}{C_4}$ containing no induced copy of $C_4$. For each $g\in V(G)$, let $W_g$ be the intersection of $W$ and $\{g\}\square H$, and let $w_g \coloneq |E(W_g)|$. For each $e=h_1h_2\in E(H)$, let $w_e$ denote the number of vertices $g\in V(G)$ such that $\{g\}\times\{h_1,h_2\}\subseteq  V(W)$. Then we have
\[\sum_{g\in V(G)}w_g=\Big\lvert\big\{\{g\}\times\{h_1,h_2\}\subseteq  V(W):\, g \in V(G),\,h_1h_2\in E(H)\big\}\Big\rvert=\sum_{e\in E(H)}w_e.\]
Note that $W_g$ is a subgraph of a copy of $H$, thus $|E(W_g)|\ge |V(W_g)|-\alpha(H)$, which implies that \[\sum_{g\in V(G)}w_g\geq\sum_{g\in V(G)}\left(|V(W_g)|-\alpha(H)\right)=|V(W)|-|V(G)|\alpha(H).\] On the other hand, we have $w_e\leq\alpha(G)$ for every $e\in E(H)$, as otherwise we can find $\{g_1,g_2\}\times\{h_1,h_2\}\subseteq  V(W)$ with $g_1g_2\in E(G)$ and $h_1h_2\in E(H)$, which yields an induced copy of $C_4$ in $W$, a contradiction. Hence, it holds that \[\sum_{e\in E(H)}w_e\leq|E(H)|\alpha(G).\] Combining the above, we obtain $|V(W)|-|V(G)|\alpha(H)\leq|E(H)|\alpha(G)$. Therefore,
\begin{equation}
\label{upper}
\alphagraphind{G\square H}{C_4}\leq|V(G)|\alpha(H)+|E(H)|\alpha(G).
\end{equation}
Plugging~\eqref{upper} in~\eqref{chromatic}, we get
\begin{equation} \label{B3}
\chigraphind{\R^n}{C_4}\geq\frac{|V(G)||V(H)|}{|V(G)|\alpha(H)+|E(H)|\alpha(G)}=\left(\frac{\alpha(H)}{|V(H)|}+\frac{|E(H)|}{|V(H)|}\cdot\frac{\alpha(G)}{|V(G)|}\right)^{-1}.
\end{equation}

We shall select the appropriate graphs $G$ and $H$ to optimize~\eqref{B3}.
We start with $G$. Recall from \Cref{sec:b1} that there is a unit-distance graph $G_n$ in $\R^n$ such that $\frac{\alpha(G_n)}{|V(G_n)|}\leq\left(1.239\ldots+o(1)\right)^{-n}$. Since this graph has the smallest known value of $\frac{\alpha(F)}{|V(F)|}$ among all unit-distance graphs $F$ in $\R^n$ (any better graph would give an improvement to~\eqref{chiRn}), we take it to be our $G$.

Now, let $H=J(n,k,t)$ be the generalized Johnson graphs defined in \Cref{sec:growing-canonical}. Note that if $uv$ is an edge of $H$, then $\|u-v\|^2=2(k-t)$, and thus, after a proper scaling of the unit cube, we observe that $H$ is a unit-distance graph in $\R^n$. It is immediate from the definition that $|V(H)|=\binom{n}{k}$ and $|E(H)|=\frac{1}{2}\binom{n}{k}\binom{k}{t}\binom{n-k}{k-t}$. Frankl and Wilson~\cite{FW81} showed that if $k \le n/2$, $t\le k/2$ and $k-t$ is a prime number, then $\alpha(H)\leq\binom{n}{k-t-1}.$
Now we take $k\approx 0.0453\cdot n, t\approx k/2$ and use Stirling's formula to approximate all the binomial coefficients. By plugging the resulting exponents back in~\eqref{B3}, we obtain the desired lower bound on $\chigraphind{\R^n}{C_4}$.
\end{proof}

\addtocontents{toc}{\protect\setcounter{tocdepth}{2}}

\end{document}